\documentclass[12pt]{article}
\usepackage{amsthm}
\usepackage{mathrsfs}
\usepackage{graphicx}
\usepackage{amsmath}
\usepackage{amsfonts}
\usepackage{amssymb}
\usepackage{url}
\usepackage{dsfont} % needed for double-struck 1
\usepackage{verbatim}
\usepackage[mathlines]{lineno}
\usepackage{color}
\definecolor{red}{rgb}{1,0,0}

\definecolor{pink}{rgb}{.9,.2,.7}

\definecolor{org}{rgb}{1,.7,0}

\definecolor{purp}{rgb}{.5,0,.5}

\definecolor{br}{rgb}{.5,.7,1}

\newtheorem{thm}{Theorem}[section]

\newtheorem{prop}[thm]{Proposition}
\newtheorem{cor}[thm]{Corollary}
\newtheorem{rem}[thm]{Remark}
\newtheorem{lem}[thm]{Lemma}

\newtheorem{obs}[thm]{Observation}

\def\noi{\noindent}
\def\mtx#1{\begin{bmatrix} #1 \end{bmatrix}}

\newcommand{\ba}{\begin{array}}
\newcommand{\ea}{\end{array}}

\newcommand{\bit}{\begin{itemize}}
\newcommand{\eit}{\end{itemize}}
\newcommand{\ben}{\begin{enumerate}}
\newcommand{\een}{\end{enumerate}}
\newcommand{\bea}{\begin{eqnarray*}}
\newcommand{\eea}{\end{eqnarray*}}
\newcommand{\beq}{\begin{equation}}
\newcommand{\eeq}{\end{equation}}
\newcommand{\bpf}{\begin{proof}}
\newcommand{\epf}{\end{proof}}

\newcommand{\epr}{\operatorname{epr}}
\newcommand{\pr}{\operatorname{pr}}

\newcommand{\rank}{\operatorname{rank}}
\DeclareMathSymbol{\mlq}{\mathord}{operators}{``}
\DeclareMathSymbol{\mrq}{\mathord}{operators}{`'}

\newcommand{\OL}{\overline}
\usepackage{enumerate}
\usepackage{enumitem}
\usepackage{hyperref}% Creates bookmarks

\usepackage{authblk} % Added to use \affil
\begin{document}

%\linenumbers

%\title{Title here\thanks{This research was supported by NSF DMS-0502354.}}
%\author{Xavier Martinez-Rivera\thanks{Department of Mathematics, Iowa State University, Ames,
%IA 50011, USA (xaviermr@iastate.edu).}
%}

\title{{\bf Classification of families of pr- and epr-sequences}}

%%%%%%%%%%%%%%Commented for LAMA%%%%%%%%%%%%%
%%%%%%%%%%%%%%Commented for LAMA%%%%%%%%%%%%%

%\author{Xavier Mart\'{i}nez-Rivera
%\thanks{Department of Mathematics, Iowa State University,
%Ames, IA 50011, USA (xaviermr@iastate.edu).}
%}

%%%%%%%%%%%%%%%%%%%%%%%%%%%%%%%%%%%%%%%%%%%%%
%%%%%%%%%%%%%%%%%%%%%%%%%%%%%%%%%%%%%%%%%%%%%

%%%%%%%%%%%%%%For LAMA%%%%%%%%%%%%%
%%%%%%%%%%%%%%For LAMA%%%%%%%%%%%%%

\author{Xavier Mart\'{i}nez-Rivera\thanks{%$^\ast$
Corresponding author.
Email: xaviermr@iastate.edu}
}

%%%%%%%%%%%%%%%%%%%%%%%%%%%%%%%%%%%
%%%%%%%%%%%%%%%%%%%%%%%%%%%%%%%%%%%

%%%%%%%%%%%%%%For LAMA%%%%%%%%%%%%%
%%%%%%%%%%%%%%For LAMA%%%%%%%%%%%%%

\affil{\textit{Department of Mathematics, Iowa State University, Ames, IA, USA}}

%%%%%%%%%%%%%%%%%%%%%%%%%%%%%%%%%%%
%%%%%%%%%%%%%%%%%%%%%%%%%%%%%%%%%%%

\maketitle

\begin{abstract}
This paper establishes new restrictions for attainable enhanced principal rank characteristic sequences (epr-sequences).
These results are then used to classify two related
families of sequences that are attainable by a real symmetric matrix:
the family of
principal rank characteristic sequences (pr-sequences)
not containing three consecutive $1$s and
the family of
epr-sequences
which contain an $\tt{N}$ in every subsequence of length $3$.
\end{abstract}

\noi{\bf Keywords:}
Principal rank characteristic sequence;
enhanced principal rank characteristic sequence;
minor; rank; symmetric matrix\\

\noi{\bf AMS Subject Classifications:}
15A15; 15A03; 15B57

%%%%%%%%%%%%%%% Introduction %%%%%%%%%%%%%%%%%%%%%

\section{Introduction}\label{S:intro}
$\indent$
Given an $n \times n$ symmetric matrix $B$ over a
field $F$,
\textit{the principal rank characteristic sequence}
(abbreviated pr-sequence) of $B$ is defined as
$\pr(B) = r_0]r_1 \cdots r_n$, where
   \begin{equation*}
      r_k =
         \begin{cases}
             1 &\text{if $B$ has a nonzero principal minor of order $k$, and}\\
             0 &\text{otherwise,}
         \end{cases}
   \end{equation*}
while $r_0 = 1$ if and only if $B$ has a
$0$ diagonal entry \cite{ORIGINAL};
the \textit{order} of a minor is $k$ if it is
the determinant of a $k \times k$ submatrix.

The \textit{principal minor assignment problem}, introduced in \cite{PMAP introduced}, asks the following question: can we find an $n \times n$ matrix with prescribed principal minors?
As a simplification of the principal minor assignment problem, Brualdi et al.\ \cite{ORIGINAL} introduced the pr-sequence of a real symmetric matrix as defined above.
An attractive result obtained in \cite{ORIGINAL} is the requirement that a pr-sequence that can be realized by a real symmetric matrix cannot contain the subsequence $001$, meaning that in the pr-sequence of such matrix,
the presence of the subsequence 00 forces 0s from that point forward.
This result was later generalized by
Barrett et al.\ \cite{FIELDS} for symmetric matrices over any field; this led them to the study of symmetric matrices over various fields, where, among other results,
a characterization of the pr-sequences that can be realized by a symmetric matrix over a field of characteristic 2 was obtained.
Although not deeply studied, the family of pr-sequences not containing three consecutive 1s were of interest in \cite{ORIGINAL}, since the pr-sequences of the principal submatrices of a matrix realizing a pr-sequence not containing three consecutive 1s
possess the rare property of being able to inherit the majority of the 1s of the original sequence;
this family will be one of the central themes of
this paper.

Due to the limitations of the pr-sequence, which only records the presence or absence of a full-rank principal submatrix of each possible order,
Butler et al.\ \cite{EPR} introduced the
the \textit{enhanced principal rank characteristic sequence} (abbreviated epr-sequence) of an $n \times n$ symmetric matrix $B$ over a field $F$, denoted by
$\epr(B) = \ell_1\ell_2 \cdots \ell_n$, where
   \begin{equation*}
      \ell_k =
         \begin{cases}
             \tt{A} &\text{if all the principal minors of order $k$ are nonzero;}\\
             \tt{S} &\text{if some but not all the principal minors of order $k$ are nonzero;}\\
             \tt{N} &\text{if none of the principal minors of order $k$ are nonzero, i.e., all are zero.}
         \end{cases}
   \end{equation*}

A (pr- or epr-) sequence is said to be
\textit{attainable} over a field $F$ provided that there exists a symmetric matrix $B \in F^{n \times n}$ that attains it;
otherwise, we say that it is \textit{unattainable}.
Among other results, techniques to construct attainable epr-sequences were presented in \cite{EPR}, as well as necessary conditions for an epr-sequence to be attainable by a symmetric matrix, with many of them asserting that subsequences such as $\tt NSA$, $\tt NAN$ and $\tt NAS$, among others, cannot occur in epr-sequences over certain fields.
Continuing the study of epr-sequences,
Fallat et al.\ \cite{skew} characterized all the epr-sequences that are attainable by skew-symmetric matrices.

In this paper the study of pr- and epr-sequences of symmetric matrices is continued.
Section \ref{s:restric} establishes new restrictions for
epr-sequences to be attainable over certain fields.
The results from Section \ref{s:restric} are then implemented in Section \ref{s:no 111}, where, for real symmetric matrices,
we classify all the attainable pr-sequences not
containing three consecutive $1$s.
Using this classification, in Section \ref{s:N in length 3}, a related family of attainable epr-sequences is classified, namely those that contain an $\tt{N}$ in every subsequence of length $3$.
We then conclude with Proposition \ref{Uniqueness: general}, where we highlight an interesting property exhibited by the vast majority of attainable
pr-sequences not containing three consecutive 1s;
that is, the property of being associated with a
unique attainable epr-sequence.

A pr-sequence and an epr-sequence are
\textit{associated} with each other if a matrix
(which may not exist) attaining the epr-sequence also
attains the pr-sequence.
A subsequence that does not appear in an attainable
sequence is \textit{forbidden}
(and we may also say that it is \textit{prohibited}).
Moreover, a sequence is said to have \textit{order} $n$ if it corresponds to a matrix of order $n$,
while a subsequence has \textit{length} $n$ if
it consists of $n$ terms.

Let $B=[b_{ij}]$ and let
$\alpha, \beta \subseteq \{1, 2, \dots, n\}$.
Then the submatrix lying in rows indexed by $\alpha$,
and columns indexed by $\beta$,
is denoted by $B[\alpha, \beta]$;
if $\alpha = \beta$, then
$B[\alpha, \alpha]$ is abbreviated to $B[\alpha]$.
The matrices $0_n$, $I_n$ and $J_n$ are the matrices
of order $n$
denoting the zero matrix, the identity matrix and
the all-$1$s matrix, respectively.
The direct sum of two matrices $B$ and $C$ is
denoted by $B \oplus C$.
Given a graph \textit{G}, $A(G)$ denotes the adjacency matrix of $G$, while $P_n$ and $C_n$ denote the path and cycle, respectively, on $n$ vertices.

\subsection{\textit{Results cited}}
$\indent$
The purpose of this section is to list results we will
cite frequently, and assign abbreviated nomenclature to
some of them.

\begin{thm}\label{Inverse pr Theorem}{\rm\cite[Theorem 2.7]{ORIGINAL}
Suppose $B$ is a nonsingular real symmetric matrix with $\pr(B) = r_0]r_1 \cdots r_n$.
Let $\pr(B^{-1}) = r'_0]r'_1 \cdots r'_n$.
Then $r'_n = r_n = 1$,
while for each $i$ with $1 \leq i \leq n-1$, $r'_i = r_{n-i}$.
Finally, $r'_0 = 1$ if and only if $B$ has some principal minor of order
$n-1$ that is zero.}
\end{thm}

\begin{thm}\label{001 is forbidden}
{\rm \cite[Theorem 4.4]{ORIGINAL}
($00$ Theorem)}
Let $B$ be a real symmetric matrix.
Let $pr(B) = r_0]r_1 \cdots r_n$ and suppose that,
for some $k$ with $0 \leq k \leq n-2$, $r_{k+1} = r_{k+2} = 0$.
Then $r_i = 0$ for all $i \geq k+1$.
In particular, $r_n = 0$, so that $B$ is singular.
\end{thm}

\begin{thm}\label{0110...1 is forbidden}
{\rm \cite[Theorem 6.5]{ORIGINAL}
($0110$ Theorem)}
Suppose $n \geq 4$ and $pr(B) = r_0]r_1 \cdots r_n$.
If, for some $k$ with $1 \leq k \leq n-3$,
$r_k = r_{k+3} = 0$, then $r_{i} = 0$ for all
$k+3 \leq i \leq n$. In particular, $B$ is singular.

\end{thm}

A generalization of Theorem \ref{001 is forbidden} in
\cite{FIELDS} led to an analogous result for
epr-sequences over any field:

\begin{thm}\label{NN result}
{\rm \cite[Theorem 2.3]{EPR}
($\tt NN$ Theorem)}
Suppose $B$ is a symmetric matrix over
a field $F$,
$\epr(B) = \ell_1\ell_2 \cdots \ell_n$, and
$\ell_k = \ell_{k+1} = \tt{N}$ for some $k$.
Then $\ell_i = \tt{N}$ for all $i \geq k$.
(That is, if an epr-sequence of a matrix ever has $\tt{NN}$, then it must have $\tt{N}$s from that point forward.)
\end{thm}

\begin{thm}\label{Inverse Thm}
{\rm \cite[Theorem 2.4]{EPR} (Inverse Theorem)}
Suppose $B$ is a nonsingular symmetric matrix
over a field $F$.
If $\epr(B) = \ell_1\ell_2 \cdots \ell_{n-1}\tt{A}$, then
$\epr(B^{-1}) = \ell_{n-1}\ell_{n-2} \cdots \ell_{1}\tt{A}$.
\end{thm}

Each instance of $\cdots$ below is permitted to be empty.

\begin{prop}\label{SN...A...}
{\rm \cite[Proposition 2.5]{EPR}}
The epr-sequence $\tt{SN}\cdots \tt{A} \cdots$ is forbidden for symmetric matrices over any field.
\end{prop}

We say that
$\tt{SN} \cdots \tt{A} \cdots$ is prohibited when referencing Proposition \ref{SN...A...}.

\begin{thm}\label{Inheritance}
{\rm{\cite[Theorem 2.6]{EPR}}}
{\rm (Inheritance Theorem)}
Suppose that  $B$ is a symmetric matrix over a field $F$, $m \leq n$, and $1\le i \le m$.  % has epr-sequence $\ell_1\ell_2\cdots\ell_n$.
\ben
\item
If $[\epr(B)]_i={\tt N}$, then  $[\epr(C)]_i={\tt N}$ for all $m\times m$ principal submatrices $C$.
\item If  $[\epr(B)]_i={\tt A}$, then  $[\epr(C)]_i={\tt A}$ for all $m\times m$ principal submatrices $C$.

\item If $[\epr(B)]_m={\tt S}$, then there exist $m\times m$ principal submatrices $C_A$ and $C_N$ of $B$ such that $[\epr(C_A)]_m = {\tt A}$ and $[\epr(C_N)]_m = {\tt N}$.
\item If $i < m$ and $[\epr(B)]_i = {\tt S}$, then there exists an $m \times m$ principal submatrix $C_S$ such that $[\epr(C_S)]_i ={\tt S}$.
\een

\end{thm}

\begin{cor}\label{NSA & ...ASN...A...}
{\rm\cite[Corollary 2.7]{EPR}}
No symmetric matrix over any field can have
$\tt{NSA}$ in its epr-sequence. Further, no symmetric matrix over any field can have the epr-sequence $\cdots \tt{ASN} \cdots \tt{A} \cdots$.
\end{cor}

Corollary \ref{NSA & ...ASN...A...} will be invoked by just stating that $\tt{NSA}$ or $\cdots \tt{ASN} \cdots \tt{A} \cdots$ is prohibited.

If $B$ is a matrix with a
nonsingular principal submatrix $B[\alpha]$,
$B/ B[\alpha]$ denotes the
Schur complement of $B[\alpha]$ in $B$.

\begin{thm} \label{schur}
{\rm{\cite[Proposition 2.13]{EPR}}}
{\rm (Schur Complement Theorem)}
Suppose $B$ is a symmetric matrix over a field of characteristic not $2$ with $\rank B=m$.
Let $B[\alpha]$ be a nonsingular
principal submatrix of $B$ with $|\alpha| = k \leq m$, and let $C = B/B[\alpha]$.
Then the following results hold.%\vspace{-3mm}
\begin{enumerate}
\item\label{p1SC} $C$ is an $(n-k)\times (n-k)$ symmetric matrix. %\vspace{-3mm}
\item\label{p2SC} Assuming the indexing of $C$ is inherited from $B$, any principal minor of $C$
is given by %\vspace{-3mm}
\[ \det C[\gamma] = \det B[\gamma \cup \alpha]/ \det B[\alpha].\]
\item\label{p3SC} $\rank C = m-k$.%\vspace{-3mm}
\item\label{LHp4SC} Any nonsingular principal submatrix of $B$ of order at most $m$ is contained in a nonsingular principal submatrix of order $m$.
\end{enumerate}
\end{thm}

\begin{thm}\label{NAN & NAS}
{\rm\cite[Theorem 2.14]{EPR}}
Neither the epr-sequences $\tt{NAN}$ nor $\tt{NAS}$ can occur as a subsequence of the epr-sequence of a symmetric matrix over a field of characteristic not $2$.
\end{thm}

We will refer to Theorem \ref{NAN & NAS} by simply stating that $\tt{NAN}$ or $\tt{NAS}$ is prohibited,
while Theorem \ref{ANS} below is referenced by stating that $\tt ANS$ `must be initial.'

\begin{thm}\label{ANS}{\rm\cite[Theorem 2.15]{EPR}}
In the epr-sequence of a symmetric matrix over a field of characteristic not $2$, the subsequence $\tt ANS$ can only occur as the initial subsequence.
\end{thm}

%%%%%%%%%%%%%%%%%%%%%%%%%%%%%%%%%%%%%%%%%%%%%%%%%%%%%%%%%%%%%
%%%%%%%%%%%%%%%%%%%%%%%%%%%%%%%%%%%%%%%%%%%%%%%%%%%%%%%%%%%%%
%%%%%%%%%%%%%%%%%%%%%%%%%%%%%%%%%%%%%%%%%%%%%%%%%%%%%%%%%%%%%
%%%%%%%%%%%%%%%%%%%%%%%%%%%%%%%%%%%%%%%%%%%%%%%%%%%%%%%%%%%%%
%%%%%%%%%%%%%%%%%%%%%%%%%%%%%%%%%%%%%%%%%%%%%%%%%%%%%%%%%%%%%
%%%%%%%%%%%%%%%%%%%%%%%%%%%%%%%%%%%%%%%%%%%%%%%%%%%%%%%%%%%%%
%%%%%%%%%%%%%%%%%%%%%%%%%%%%%%%%%%%%%%%%%%%%%%%%%%%%%%%%%%%%%

\section{Restrictions on attainable epr-sequences}\label{s:restric}
$\indent$
In this section, we establish new restrictions on
attainable epr-sequences.
We begin with restrictions that apply to fields of characteristic not $2$.
For convenience, given a matrix $B$,
we adopt some of the notation in \cite{ORIGINAL},
and denote with
$B_{i_1 i_2 \dots  i_k}$, the principal minor
$\det(B[\{ i_1, i_2, \dots, i_k \}])$.

\begin{prop}\label{NSNA...}
Let $n \geq 6$.
Then no $n \times n$ symmetric matrix over
a field of characteristic not $2$ has an
epr-sequence starting $\tt{NSNA} \cdots$.
\end{prop}

\bpf
Let $B = [b_{ij}]$ be an $n\times n$ symmetric matrix over a field of characteristic not $2$ and
let $\epr(B) = \ell_1\ell_2 \cdots \ell_n$.
Suppose to the contrary that $\epr(B) = \tt{NSNA}\cdots$.
Since $\ell_3 = \tt{N}$, and because
$B_{pqr} = 2b_{pq}b_{pr}b_{qr}$ for any distinct
$p,q,r \in \{1, 2,\dots, n\}$, $B[\{1,2,3\}]$ and $B[\{4,5,6\}]$
must each contain a zero off-diagonal entry.
Moreover, since $\ell_4 = \tt{A}$, $0_3$ is not
a principal submatrix of $B$,
implying that $B[\{1,2,3\}]$ and $B[\{4,5,6\}]$
must each contain a nonzero off-diagonal entry.
Since $\{1,2,3\}$ and $\{4,5,6\}$ are disjoint,
and because a simultaneous permutation of the rows and columns of a matrix has no effect on its determinant,
we may assume, without loss of generality, that
$b_{12} = b_{56} = 0$ and that
$b_{13}, b_{46}$ are nonzero.
Similarly, since $\{1,2,3\}$ and $\{4,5,6\}$ are disjoint, and because multiplication of any row and column of a matrix by a nonzero constant preserves the rank of every submatrix,
we may also assume, without loss of generality, that
$b_{13} = b_{46} = 1$.
We consider two cases.

\textit{Case 1}: $b_{14} = 0$.
Since $\ell_4 = \tt{A}$,
$(b_{15}b_{24})^2 = B_{1245}\neq 0$;
it follows that $b_{15}$ and $b_{24}$ are nonzero.
Since $\ell_3 = \tt{N}$,
$B_{135} = 2b_{15}b_{35} = 0$;
hence, $b_{35} = 0$.
Since $B[\{3, 5, 6\}] \neq 0_3$, $b_{36} \neq 0$.
Since $ 2b_{16}b_{36} = B_{136} = 0$,
$b_{16} = 0$.
Then, as $B[\{1, 2, 6\}] \neq 0_3$, $b_{26} \neq 0$.
It follows that
$B_{246} = 2b_{24}b_{26} \neq 0$,
a contradiction to $\ell_3 = \tt{N}$,
implying that it is impossible to have $b_{14} = 0$.

\textit{Case 2}: $b_{14} \neq 0$.
Since $2b_{14}b_{34} = B_{134} = 0$, and because
$2b_{14}b_{16} = B_{146} = 0$,
$b_{34} = b_{16} = 0$.
Since $B[\{1, 2, 6\}] \neq 0_3$, $b_{26} \neq 0$.
Since $ 2b_{24}b_{26} = B_{246} = 0$,
$b_{24} = 0$.
Since $(b_{14}b_{23})^2 = B_{1234} \neq 0$,
$b_{23} \neq 0$.
Then, as $ 2b_{23}b_{26}b_{36} = B_{236} = 0$,
$b_{36} = 0$.
It follows that $B_{1356} = 0$,
a contradiction to $\ell_4 = \tt{A}$.
\epf

It should be noted that
$\tt{NSNA}$ and $\tt{NSNAA}$ are attainable by
$A(P_4)$ and $A(C_5)$, respectively \cite{EPR},
but this does not contradict Proposition \ref{NSNA...},
which requires $n \geq 6$.

\begin{prop}\label{NSNA in the first n-2 terms}
Let $B$ be a symmetric matrix over
a field of characteristic not $2$ and
$\epr(B) = \ell_1 \ell_2 \cdots \ell_n$.
Then $\tt{NSNA}$ cannot occur as a subsequence of
$\ell_1 \ell_2 \cdots \ell_{n-2}$.
\end{prop}

\bpf
If $n \leq 5$, the result follows vacuously.
So, assume $n \geq 6$.
Suppose to the contrary that $\tt{NSNA}$ occurs as a subsequence of $\ell_1\ell_2 \cdots \ell_{n-2}$ and that
$\ell_{k} \ell_{k+1} \ell_{k+2} \ell_{k+3} = \tt{NSNA}$,
for some $k$ with $1 \leq k \leq n-5$.
By Proposition \ref{NSNA...}, $k \geq 2$, and,
by the $\tt{NN}$ Theorem, $\ell_{k-1} \neq \tt{N}$;
it follows that $B$ has a $(k-1) \times (k-1)$
nonsingular principal submatrix, say $B[\alpha]$.
By the Schur Complement Theorem, $B/B[\alpha]$ has an
epr-sequence starting $\tt{NXNAYZ} \cdots$,
where $\tt{X}, \tt{Y}, \tt{Z} \in \{A,S,N\}$.
The $\tt{NN}$ Theorem and the fact that
$\tt{NAN}$ is prohibited imply that
$\tt{X} = \tt{S}$; hence,
$\epr(B)$ starts $\tt{NSNAYZ} \cdots$,
a contradiction to Proposition \ref{NSNA...}.
\epf

With the next result,
we generalize (and provide a simpler proof of) \cite[Proposition 2.11]{EPR}.

\begin{prop}\label{SAN}Suppose $B$ is a symmetric matrix
over a field of characteristic not $2$,
$\epr(B) = \ell_1 \ell_2 \cdots \ell_n$ and
$\ell_k \ell_{k+1}\ell_{k+2} = \tt{SAN}$ for some $k$.
Then $\ell_j = \tt{N}$ for all $j \geq k+2$.
\end{prop}

\bpf
If $n=3$, we are done. Suppose $n>3$.
Suppose that $\ell_k\ell_{k+1}\ell_{k+2} = \tt{SAN}$
for some $k$ with $1 \leq k \leq n-2$.
If $k=n-2$, we are done. Suppose $k < n-2$.
By \cite[Corollary 2.10]{EPR},
which prohibits $\tt{SANA}$,
$\ell_{k+3} \neq \tt{A}$.
Since $\tt{ANS}$ must be initial, $\ell_{k+3} \neq \tt{S}$. Hence, $\ell_{k+3} = \tt{N}$.
The desired conclusion now follows from
the $\tt{NN}$ Theorem.
\epf

We now confine our attention to real symmetric matrices.
The next result is immediate from
Theorem \ref{0110...1 is forbidden}.

\begin{prop}\label{NXXN}
Let $B$ be a real symmetric matrix and
$\epr(B) = \ell_1\ell_2 \cdots \ell_n$.
Suppose $\ell_k = \ell_{k+3}=\tt{N}$ for some $k \geq 1$.
Then $\ell_i = \tt{N}$ for all $i \geq k+3$.
In particular, $B$ is singular.
\end{prop}

We emphasize that Proposition \ref{NXXN}
asserts that a sequence of the form
$\cdots \tt{N}\tt{X}\tt{Y}\tt{N} \cdots \tt{Z} \cdots$,
with $\tt{X},\tt{Y} \in \{ \tt{A}, \tt{S}, \tt{N}\}$ and
$\tt{Z} \in \{ \tt{A},\tt{S}\}$, is not attainable by a real symmetric matrix.

%\begin{cor}\label{N in every subseq of order 2}
%Let $n \geq 4$, $B$ be a real symmetric matrix, and
%$\epr(B) = \ell_1\ell_2 \cdots \ell_n$.
%Suppose $\epr(B)$ contains an $\tt{N}$ in
%every subsequence of length $3$.
%Then, for all $k$ such that $2 \leq k \leq \rank(B)-2$, $\ell_k\ell_{k+1}$
%contains an $\tt{N}$.
%\end{cor}

Given a sequence
$t_{i_{1}} t_{i_{2}} \cdots t_{i_{k}}$,
$\overline{t_{i_{1}} t_{i_{2}} \cdots t_{i_{k}}}$
indicates that the sequence may be repeated as many
times as desired (or it may be omitted entirely).
According to \cite[Proposition 2.17]{EPR},
the sequence $\tt{ANA\overline{A}}$ is attainable by a symmetric matrix over a field of characteristic $0$. \cite[Table 1]{EPR} raises the following question:
does a real symmetric matrix,
with an epr-sequence starting $\tt{ANA \cdots}$,
always have epr-sequence $\tt{ANA\overline{A}}$? The answer is affirmative;
what follows makes this precise.

\begin{prop}\label{ANA}
Any $n \times n$ real symmetric matrix with an epr-sequence starting $\tt{ANA \cdots}$ is conjugate by a nonsingular
diagonal matrix to one of $\pm(J_n-2I_n)$.
Furthermore, its epr-sequence is $\tt{ANA\overline{A}}$.
\end{prop}

\bpf Let $B = [b_{ij}]$ be an
$n \times n$ real symmetric matrix with an
epr-sequence starting $\tt{ANA \cdots}$.
Notice that all the diagonal entries of $B$ must have
the same sign, as otherwise there would be a
principal minor of order $2$ that is nonzero.
Let $C = [c_{ij}]$ be the matrix among $B$ and $-B$ with all diagonal entries negative.
Let $D = [d_{ij}]$ be the $n \times n$ diagonal matrix
with $d_{11} = 1/\sqrt{-c_{11}}$ and
$d_{jj} = \rm{sign}(c_{1j})/\sqrt{-c_{jj}}$ for $j \geq 2$.
Now, notice that every entry of $DCD$ is $\pm 1$,
every diagonal entry is $-1$ and
every off-diagonal entry in the first row and
the first column is $1$.
We now show that $DCD = J_n - 2I_n$.
Since multiplication of any row and column of a matrix
by a nonzero constant preserves the rank of every submatrix,
$\epr(DCD) = \epr(C) = \epr(B)$.
Let $i,j \in \{2,3, \dots n\}$ be distinct,
$\alpha = \{1, i, j \}$ and let
$a$ be the $(i,j)$-entry of $DCD$.
A simple computation shows that
$\det((DCD)[\alpha]) = (a+1)^2$.
Since every principal minor of order $3$ of
$DCD$ is nonzero, $a = 1$.
Then, as $i$ and $j$ were arbitrary,
$DCD = J_n - 2I_n$.
Then, as $C=B$ or $C=-B$,
it follows that $B$ is conjugate by a nonsingular
diagonal matrix to one of $\pm(J_n-2I_n)$, and that
$\epr(B) = \epr(J_n - 2I_n) = \tt{ANA\overline{A}}$
(see \cite[Proposition 2.17]{EPR}).
\epf

We are now in position to prove the following result.

\begin{thm}\label{(A)ANAA(A)}
Any epr-sequence of a real symmetric matrix containing
$\tt{ANA}$ as a non-terminal subsequence is of the form $\tt{\overline{A}ANAA\overline{A}}$.
\end{thm}

\bpf Let $B$ be a real symmetric matrix containing
$\tt{ANA}$ as a non-terminal subsequence.
Let $\epr(B) = \ell_1\ell_2 \cdots \ell_n$.
Suppose $\ell_{k+1}\ell_{k+2}\ell_{k+3} = \tt{ANA}$
for some $k$ with $0 \leq k \leq n-4$.
Since $\tt{NAN}$ and $\tt{NAS}$ are prohibited,
$\ell_{k+4} = \tt{A}$.
If $k=0$, the conclusion follows from Proposition \ref{ANA};
so, assume $k>0$.
Suppose $\ell_i \neq \tt{A}$ for some $i$ with $i < k+1$.
By the Inheritance Theorem,
$B$ has a (nonsingular) $(k+4) \times (k+4)$ principal submatrix $B'$ whose epr-sequence
$\ell'_1\ell'_2 \cdots \ell'_{k+4}$ ends with $\tt{ANAA}$
and has $\ell'_i \neq \tt A$.
Then, by the Inverse Theorem,
$\epr((B')^{-1})$ starts with $\tt{ANA}$ and
$\epr((B')^{-1}) \neq \tt{ANA\overline{A}}$,
a contradiction to Proposition \ref{ANA}.
Thus, $\epr(B) = \tt{\overline{A}AANAA}$$\ell_{k+5} \cdots \ell_n$,
where $\ell_{k+5} \cdots \ell_n$ may not exist.

We now show that
$\ell_{k+5} \cdots \ell_n = \tt \OL{A}$.
If $n = k+4$, we are done; so, suppose $n > k+4$.
We proceed by contradiction, and consider two cases.

\textit{Case 1}: $\ell_j = \tt{N}$ for some $j > k+4$.
Since $\ell_k = \tt{A}$, there exists a $k \times k$ principal submatrix of $B$, say $B[\alpha]$, that is nonsingular.
Let $C = B/B[\alpha]$.
By the Schur Complement Theorem,
$C$ has order $n-k$, $\epr(C)$ starts $\tt{ANA \cdots}$ and $\epr(C)$ has an $\tt{N}$ in the ($j-k$)-th position;
hence, $\epr(C) \neq \tt{ANA\overline{A}}$,
a contradiction to Proposition \ref{ANA}.
It follows that a sequence containing $\tt{ANA}$ as a
non-terminal subsequence cannot contain an $\tt{N}$ from that point forward, implying that any real symmetric matrix with an epr-sequence containing
$\tt{ANA}$ is nonsingular.

\textit{Case 2}: $\ell_j = \tt{S}$ for some $j > k+4$.
By the Inheritance Theorem, $B$ has a singular $j \times j$ principal submatrix whose epr-sequence contains $\tt{ANA}$, which contradicts the assertion above.

We conclude that we must have
$\ell_{k+5} \cdots \ell_n = \tt \OL{A}$, which completes the proof.
\epf

It is natural to now ask,
does Theorem \ref{(A)ANAA(A)} hold if $\tt{ANA}$ occurs at
the end of the sequence?
According to \cite[Table 5]{EPR}, $\tt{SAANA}$ is attainable,
answering the question negatively.
%(SAANA that is the only counterexample among sequences of order $n \leq 5$).

\begin{thm}\label{SNA}
Let $B$ be a real symmetric matrix with
$\epr(B) = \ell_1\ell_2 \cdots \ell_n$.
Then $\tt{SNA}$ cannot occur as a subsequence of
$\ell_1\ell_2 \cdots \ell_{n-2}$.
\end{thm}

\bpf
If $n \leq 4$, the result follows vacuously.
So, assume $n > 4$.
Suppose to the contrary that $\tt{SNA}$ occurs as a subsequence of $\ell_1\ell_2 \cdots \ell_{n-2}$, and that
$\ell_{k+1}\ell_{k+2}\ell_{k+3} = \tt{SNA}$ for
some $k$ with $0 \leq k \leq n-5$.
Since  $\tt{SN} \cdots \tt{A} \cdots$ is prohibited,
$k \geq 1$.
Since $\tt{NAN}$ and $\tt{NAS}$ are prohibited,
$\ell_{k+4} = \tt{A}$.
Then, as $\tt{ASNA}$ is prohibited,
$\ell_k \neq \tt{A}$.
And, by Proposition \ref{NSNA in the first n-2 terms},
$\ell_k \neq \tt N$;
it follows that $\ell_k = \tt S$.
Thus, we have $\ell_{k} \cdots \ell_{k+4} = \tt{SSNAA}$.
We examine the three possibilities for $\ell_{k+5}$ .

\textit{Case 1}: $\ell_{k+5} = \tt{A}$.
Now we have $\ell_{k} \cdots \ell_{k+5} = \tt{SSNAAA}$.
By the Inheritance Theorem,
$B$ has a $(k+5) \times (k+5)$ principal submatrix $B'$
whose epr-sequence ends with $\tt{SXNAAA}$, where
$\tt{X} \in \{\tt{A}, \tt{S}, \tt{N}\}$.
By the $\tt{NN}$ Theorem, $\tt{X} \neq \tt{N}$; and,
by Proposition \ref{SAN}, $\tt{X} \neq \tt{A}$;
it follows that $\tt{X=S}$.
By the Inverse Theorem, $\epr((B')^{-1})$
contains $\tt{ANS}$ as a non-initial subsequence,
a contradiction, since $\tt{ANS}$ must be initial.
We conclude that $\ell_{k+5} \neq \tt{A}$.

\textit{Case 2}: $\ell_{k+5} = \tt{N}$.
Now we have $\ell_{k} \cdots \ell_{k+5} = \tt{SSNAAN}$.
Since $\ell_k = \tt{S}$,
$B$ has a $k \times k$ nonsingular principal submatrix, say $B[\alpha]$.
By the Schur Complement Theorem, $B/B[\alpha]$ has an epr-sequence starting $\tt{YNAAN} \cdots$,
where $\tt{Y} \in \{\tt{A, S, N}\}$.
By Theorem \ref{(A)ANAA(A)}, $\tt{Y} \neq \tt{A}$;
since $\tt{SN} \cdots \tt{A} \cdots$ is prohibited, $\tt{Y} \neq \tt{S}$; and,
by the $\tt{NN}$ Theorem, $\tt{Y} \neq \tt{N}$.
It follows that we must have
$\ell_{k+5} \neq \tt{N}$.

From Cases 1 and 2 we can deduce that
the subsequence $\tt{SSNAAZ}$, where
$\tt{Z} \in \{\tt{A, N}\}$, cannot occur in
the epr-sequence of a real symmetric matrix.

\textit{Case 3}: $\ell_{k+5} = \tt{S}$.
Now we have $\ell_{k} \cdots \ell_{k+5} = \tt{SSNAAS}$.
By the Inheritance Theorem,
$B$ has a $(k+5) \times (k+5)$
principal submatrix with an epr-sequence ending with
$\tt{SXNAAY}$, where $\tt{X} \in \{\tt{A, S, N}\}$ and
$\tt{Y} \in \{\tt{A, N}\}$.
By the $\tt{NN}$ Theorem, $\tt{X} \neq \tt{N}$;
and, by Proposition \ref{SAN}, $\tt{X} \neq \tt{A}$.
It follows that $\tt{X=S}$,
which contradicts the assertion above.
\epf

As $\tt{NAN}$ is prohibited, we have the following corollary to Theorem \ref{SNA}.

\begin{cor}\label{SNA and SNAA} The only way $\tt{SNA}$ can occur in
the epr-sequence of a real symmetric matrix is in one of
the two terminal sequences $\tt{SNA}$ or $\tt{SNAA}$.
\end{cor}

We note that the epr-sequences $\tt{ANSSSNA}$ and $\tt{ANSSSNAA}$ are attainable \cite[Table 1]{EPR},
implying that $\tt{SNA}$ is not completely prohibited in the epr-sequence of a real symmetric matrix.
Theorem \ref{(A)ANAA(A)} and Corollary \ref{SNA and SNAA} lead to the following observation.

\begin{obs}
Any epr-sequence of a real symmetric matrix that contains $\tt{NA}$ as a non-initial subsequence is of the form $\tt{\cdots NA \overline{A}}$.
\end{obs}

The following results in this section will be of particular relevance to the main results in
Sections \ref{s:no 111} and \ref{s:N in length 3}.

\begin{lem}\label{J n/2 +1}% B doesn't have to be symmetric
Let $n$ be even and
$B$ be a nonsingular $n \times n$ real symmetric matrix.
Then $J_{\frac{n}{2} +1}$ is not a principal submatrix of B.
\end{lem}

\bpf
For the sake of contradiction, suppose, without loss of generality, that $B[\{1, \dots, \frac{n}{2}+1\}] = J_{\frac{n}{2} +1}$. Then the rank of the matrix consisting of the first $\frac{n}{2} +1$ columns of $B$ has rank less than $\frac{n}{2} +1$;
hence, $B$ is singular, a contradiction.
\epf

\begin{lem}\label{n/2 + 1pos entries}
Let $n \geq 8$ be even.
Let $B$ be an $n \times n$ nonsingular real symmetric matrix with every entry $\pm 1$ and
all entries in the first row, the first column, and
the diagonal equal to $1$.
Suppose that $\epr(B) = \ell_1\ell_2 \cdots \ell_n$ and
that  $\ell_4 = \tt{N}$.
Then every row and column of $B$ has
at most $\frac{n}{2} -1$ negative entries.
Equivalently,
every row and column of $B$ has at least
$\frac{n}{2} + 1$ positive entries.
\end{lem}

\bpf
Suppose $B = [b_{ij}]$ contains a row with
$\frac{n}{2}$ negative entries.
Let $U = \{3,4, \dots, \frac{n}{2} + 2\}$.
Without loss of generality, suppose
$b_{2j} = -1$ for all $j \in U$.
We claim that $B[\{1\} \cup U] = J_{\frac{n}{2} +1}$.
Suppose to the contrary that $B[\{1\} \cup U] $ contains a
negative entry;
without loss of generality, we may assume that
this entry is $b_{34}$. It follows that $B[\{1,2,3,4\}]$ is
nonsingular, a contradiction to $\ell_4 = \tt{N}$;
hence, $B[\{1\} \cup U] = J_{\frac{n}{2} +1}$.
By Lemma \ref{J n/2 +1}, $B$ is singular,
a contradiction to the nonsingularity of $B$.
We conclude that every row and column of $B$ has at most $\frac{n}{2} -1$ negative entries.
\epf

\begin{thm}\label{ANSNSN...A}
Let $n \geq 8$ be even and $B$ be an
$n \times n$ real symmetric matrix.
Suppose that $\epr(B) = \tt{ANSNSN \cdots }$.
Then $B$ is singular.
\end{thm}

\bpf
Suppose to the contrary that $B$ is nonsingular.
Let $B = [b_{ij}]$.
By \cite[Proposition 8.1]{ORIGINAL}, we may assume that
every entry of $B$ is $\pm 1$ and
all entries in the first row, the first column, and
the diagonal are equal to $1$.
By Lemma \ref{n/2 + 1pos entries},
every row and column of $B$ has at least
$\frac{n}{2} + 1$ positive entries.
Because a simultaneous permutation of the rows and
columns of a matrix has no effect on its determinant,
we may assume, without loss of generality, that
the first $\frac{n}{2} + 1$ entries in the second row
(and column) are positive.
Let
\[\begin{array}{cc}
M_1=\mtx{
     1 &  1 &   1 &  1 \\
     1 &  1 &  -1 & -1 \\
     1 & -1 &   1 & -1 \\
     1 & -1 &  -1 &  1}\;\;\;\;\rm{and}\;\;&
M_2=\mtx{
     1 &  1 &   1 & -1 \\
     1 &  1 &  -1 &  1 \\
     1 & -1 &   1 &  1 \\
    -1 &  1 &   1 &  1}.
\end{array}\]
Since $M_1$ and $M_2$ are nonsingular, they are not
principal submatrices of $B$.
We now show by induction on the number of negative
entries in the second row that $B$ contains a
row with $\frac{n}{2}$ negative entries.
For the base case, first notice that the nonsingularity of
$B$ implies that $B$ must have a row with at least one
negative entry, as otherwise it will have a repeated row;
without loss of generality, we assume that $b_{2n} = -1$.
By Lemma \ref{J n/2 +1}, $B[\{1, \dots , \frac{n}{2} +1 \}]$
has a negative entry; without loss of generality, suppose
$b_{34} = -1$.
Then, as
$B[\{2,3,4,n \}] \neq M_2$, either
$b_{3n}$ or $b_{4n}$ is negative, implying that either
the third or fourth row contains two negative entries.
It follows that $B$ must contain a row with two
negative entries.

Now, for the inductive step, suppose
the second row contains
$2 \leq k  \leq \frac{n}{2} - 1$ negative entries.
Without loss of generality, suppose
$b_{2j} = -1$ for $j \in U = \{n-k+1, \dots , n \}$.
As in the base case, Lemma \ref{J n/2 +1} implies that
$B[\{1, \dots , \frac{n}{2} +1 \}]$ has a negative entry,
and, again, without loss of generality,
we may assume that
$b_{34} = -1$.
Since $B[\{1,2,p,q \}] \neq M_1$ for $p,q \in U$,
$b_{pq} = 1$ for all $p,q \in U$.
Similarly,
$B[\{1,3,4,j \}] \neq M_1$ and
$B[\{2,3,4,j \}] \neq M_2$ for $j \in U$, implying that
$b_{3j} \neq b_{4j}$ for all $j \in U$;
so, suppose $b_{3j} = x_j$ and $b_{4j} = -x_j$ for all
$j\in U$.
Then, as $\ell_6 = \tt{N}$, $-16(x_p - x_q)^2 =
\det B[\{1,2,3,4,p,q\}] = 0$ for all $p,q \in U$;
hence, $x_p = x_q$ for all $p,q \in U$.
It follows that either the third or the fourth row
contains $(n - (n-k+1) +1) + 1 =  k+1$ negative entries.
Hence, by induction, $B$ most have a row with
$\frac{n}{2}$ negative entries;
by Lemma \ref{n/2 + 1pos entries}, $B$ is singular,
a contradiction.
\epf

We note that Theorem \ref{ANSNSN...A} cannot be
generalized for $n$ odd, since, by the Inverse Theorem,
$\tt{AN \overline{SN}A}$ is attained by
$(A(C_n))^{-1}$ (see \cite[Observation 3.1]{EPR}).

\begin{prop}\label{SSNSNSS}
No real symmetric matrix
has an epr-sequence starting $\tt{SSNSNSS} \cdots$.
\end{prop}

\bpf
Let $B = [b_{ij}]$ be a real symmetric
with an epr-sequence starting $\tt{SSNSNSS} \cdots$.
By the Inheritance Theorem, $B$ has a  $7 \times 7$
principal submatrix $B[\alpha]$ with epr-sequence
$\ell'_1 \ell'_2$$\tt{N}$$\ell'_4$$\tt{N}$$\ell'_6 \tt{A}$.
Without loss of generality, suppose
$\alpha = \{2,3, \dots , 8 \}$.
By the $\tt{NN}$ Theorem,
$\ell'_2$, $\ell'_4$, $\ell'_6$ are not $\tt{N}$.
Since $\tt{NAN}$ and $\tt{NSA}$ are prohibited,
$\ell'_4 = \tt{S}$ and $\ell'_6 = \tt{A}$.
Since $\tt{ANS}$ must be initial, $\ell'_2 = \tt{S}$.
Hence, $\epr(B[\alpha]) = \ell'_1 \tt{SNSNAA}$.
Since $\tt{ASN} \cdots \tt{A}$ is prohibited,
$\ell'_1 \neq \tt{A}$.
Then, as the epr-sequence $\tt{SSNSNAA}$ is associated with
the pr-sequence $1]1101011$, which is unattainable by \cite[Proposition 4.1]{FIELDS},
$\tt{\ell'_1} \neq \tt{S}$;
hence, $\ell'_1 = \tt{N}$, so that
$\epr(B[\alpha]) = \tt{NSNSNAA}$.
We note that a simultaneous permutation of the rows and
columns of a matrix has no effect on its determinant;
thus, since all diagonal entries of $B[\alpha]$ are zero,
and because $B$ contains a nonzero diagonal entry,
we may assume, without loss of generality, that
$b_{11} \neq 0$.

Let $C = B[\{1 \} \cup \alpha]$ and $C=[c_{ij}]$.
Then $\epr(C)$ starts with $\tt{S}$ and
$\epr(C[\alpha]) = \epr(B[\alpha]) = \tt{NSNSNAA}$.
%%%by the Inheritance Theorem,
%%%$\epr(C) =
%%%\tt{S} \ell''_2 \tt{N} \ell''_4 \tt{N}
%%%\ell''_6 \ell''_7 \ell''_8$.
%%%Since $B[\alpha]$ is a principal submatrix of $C$,
%%%$\epr(C) =\tt{SSNSN\ell''_6 \ell''_7 \ell''_8}$ and
%%%$\tt{\ell''_6}$, $\tt{\ell''_7}$ $\in \{ \tt{A}, \tt{S}\}$.
%%%By Theorem \ref{SNA}, $\tt{\ell''_6} = \tt{S}$, and,
%%%since $\tt{NSA}$ is prohibited,
%%%$\tt{\ell''_7} = \tt{S}$.
%%%Hence, $\epr(C) = \tt{SSNSNSS\ell''_8}$.\\
Since every $6 \times 6$ principal submatrix of
$C[\alpha]$ is nonsingular,
$C[\alpha]$ contains at least two nonzero entries in
each row (and column), as otherwise $C[\alpha]$ contains a
$6 \times 6$ principal submatrix with a row (and column)
consisting of only zeros.
Moreover, we note that $c_{11} = b_{11} \neq 0$;
because multiplication of any row and column of a matrix by a nonzero constant preserves the rank of every submatrix,
we may assume without loss of generality that $c_{11} = 1$.
Since $C[\alpha]$ contains a nonzero principal minor of order $2$,
we may assume, without loss of generality, that $\det((C[\alpha])[\{1,2\}]) \neq 0$;
thus, $-(c_{23})^2 = C_{23} = (C[\alpha])[\{1,2\}] \neq 0$;
hence, $c_{23} \neq 0$, and,
without loss of generality, we may assume that
$c_{23} = 1$.
Since $C[\alpha]$ contains at least two nonzero
entries in each row and column, $c_{2j} \neq 0$ for some
$j \in \{4,5,6,7,8 \}$;
so, we may assume that $c_{24} = 1$.
It follows that $2c_{34} = C_{234} = 0$, and
so $c_{34} = 0$.
Then, as $C[\alpha]$ contains at least two nonzero
entries in each row and column, $c_{3j} \neq 0$ for some
$j \in \{5,6,7,8 \}$;
thus,
suppose $c_{35} = 1$.
It follows that $2c_{25} = C_{235} = 0$, and
so $c_{25} = 0$.
Now we have $-1 + 2c_{12}c_{13} = C_{123} = 0$,
$-1 + 2c_{12}c_{14} = C_{124} = 0$ and
$-1 + 2c_{13}c_{15} = C_{135} = 0$;
it follows that
$c_{12}$, $c_{13}$, $c_{14}$ and $c_{15}$ are nonzero.
Let $c_{12} = x$; then
$c_{13} = c_{14} = 1/2x$ and
$c_{15} = x$.
We now show that each of
$c_{16}$, $c_{17}$ and $c_{18}$ is nonzero.
Suppose to the contrary that $c_{1j} = 0$ for
some $j \in \{6,7,8\}$;
then $-(c_{ij})^2 = C_{1ij} = 0$ for all
$i \in \{3, 4, \dots, 8 \}\setminus \{j\}$;
hence, $c_{ij} = 0$ for all
$i \in \{3, 4, \dots, 8 \}$,
implying that $C[\alpha]$ contains a row with only one
nonzero entry, which is a contradiction.
Without loss of generality, we may assume that
$c_{16} = c_{17} = c_{18} = 1$.
Now, observe that $C_{145} = c_{45}(1 - c_{45})$;
since all the principal minors of order $3$ are zero, it follows that
$c_{45} =0$ or $c_{45} = 1$.
Besides for the $(1,2)$-entry $x$,
we have similar restrictions for all
the remaining unknown entries of $C$;
notice that, for $j \in \{6,7,8\}$,
$C_{12j} = c_{2j}(2x - c_{2j})$,
$C_{13j} = c_{3j}(1/x - c_{3j})$,
$C_{14j} = c_{4j}(1/x - c_{4j})$ and
$C_{15j} = c_{5j}(2x - c_{5j})$.
Similarly, for $k \in \{7,8\}$,
$C_{16k} = c_{6k}(2 - c_{6k})$.
Lastly, $C_{178} = c_{78}(2 - c_{78})$.
It is now clear that, besides the $(1,2)$-entry $x$,
each unknown entry of $C$ is restricted to
exactly two values.

We now show that $c_{45} = 1$.
Suppose to the contrary that $c_{45} = 0$.
Since $C[\alpha]$ must contain at least two nonzero entries in each row and column, without loss of generality,
we may assume that $b_{56}$ is nonzero,
implying that $c_{56} = 2x$.
Then $4xc_{36} = C_{356} = 0$, and therefore
$c_{36} = 0$.
We proceed by examining the only two possibilities for
the entry $c_{26}$. First, suppose $c_{26} = 0$.
Since all the principal minors of order $5$ of $C$ are zero,
$4xc_{46} = C_{23456} = 0$, implying that $c_{46} = 0$.
Then $C_{12456} = -4x^2 \neq 0$, a contradiction.
So, suppose $c_{26} = 2x$.
Since $4xc_{46} = C_{246} = 0$, $c_{46} = 0$.
Since $C[\alpha]$ must contain at least two nonzero entries
in each row and column, suppose, without loss of generality,
that $c_{47} \neq 0$;
hence, $c_{47} = 1/x$.
Since $2c_{27}/x = C_{247} = 0$, $c_{27}=0$.
Now, observe that
$C_{13457}=
(-2x + 2x^2c_{37} + c_{57} - xc_{37}c_{57})/2x^3$ and
$C_{23457} =
2c_{57}/x - 2c_{37}c_{57}$;
since $C_{13457}= 0$, at least one of
$c_{37}$ and $c_{57}$ is nonzero;
then, as $C_{23457}=0$, $c_{37} \neq 0$,
and so $c_{37} = 1/x$.
It follows that
$2c_{57}/x = C_{357} = 0$, and so $c_{57} = 0$.
As $-4 + 2c_{67} = C_{14567} = 0$, $c_{67} = 2$.
Then we have $C_{234567} = 0$, implying that
$C[\alpha]$ has a singular $6 \times 6$
principal submatrix, which is a contradiction.
We conclude that $c_{45} \neq 0$;
hence, $c_{45} = 1$.

Now, observe that at least one of
$c_{36}$, $c_{37}$, $c_{38}$,
$c_{46}$, $c_{47}$ and $c_{48}$ is nonzero,
as otherwise $C[{\alpha}]$, which is nonsingular,
would have two identical rows;
thus, without loss of generality, we assume that
$c_{36} \neq 0$;
hence, $c_{36} = 1/x$.
Similarly, at least one of
$c_{27}$, $c_{28}$, $c_{57}$ and $c_{58}$ is nonzero,
as otherwise
$C[\{2,3,4,5,7,8\}]=(C[\alpha])[\{1,2,3,4,6,7\}]$,
which is nonsingular, would have two identical rows;
without loss of generality,
we assume that $c_{27} \neq 0$;
thus, $c_{27} = 2x$.
Now the conditions
$C_{236} = C_{237} = C_{247} = C_{356} = 0$ imply that
$c_{26} = c_{37} = c_{47} = c_{56} = 0$.

Finally, we consider the only two possibilities for
the entry $c_{57}$. First, suppose $c_{57} = 2x$.
Then $C_{234567} = 0$, a contradiction.
Now, suppose $c_{57}= 0$.
Since $C_{234567} = -4x^2(c_{46} - 1/x)^2$ is nonzero,
$c_{46} = 0$.
Then $-2c_{67} = C_{14567}= 0$, and so $c_{67} = 0$.
Since every row and column of $C[\alpha]$ must
contain at least two nonzero entries,
it follows that $c_{68}$ and $c_{78}$ are nonzero,
implying that $c_{68} = c_{78} = 2$.
The conditions $C_{278} = C_{368} = 0$ imply that
$c_{28} = c_{38} = 0$.
Hence, $C_{23678} = 16 \neq 0$, a contradiction.
\epf

%%%%%%%%%%%%%%%%%%%%%%%%%%%%%%%%%%%%%%%%%%%%%%%%%%%%%%%%%%%%%
%%%%%%%%%%%%%%%%%%%%%%%%%%%%%%%%%%%%%%%%%%%%%%%%%%%%%%%%%%%%%
%%%%%%%%%%%%%%%%%%%%%%%%%%%%%%%%%%%%%%%%%%%%%%%%%%%%%%%%%%%%%
%%%%%%%%%%%%%%%%%%%%%%%%%%%%%%%%%%%%%%%%%%%%%%%%%%%%%%%%%%%%%
%%%%%%%%%%%%%%%%%%%%%%%%%%%%%%%%%%%%%%%%%%%%%%%%%%%%%%%%%%%%%
%%%%%%%%%%%%%%%%%%%%%%%%%%%%%%%%%%%%%%%%%%%%%%%%%%%%%%%%%%%%%

\section{Pr-sequences not containing three consecutive 1s}\label{s:no 111}
$\indent$
We begin with results that forbid certain pr-sequences
not containing three consecutive 1s;
we then implement these in
Theorem \ref{No-111 characterization}, where,
for real symmetric matrices, we classify all
the attainable pr-sequences not containing
three consecutive 1s.

It is obvious from Theorem \ref{Inverse pr Theorem} that,
with the exception of the 0th term $r'_0$,
we can explicitly determine each term in
the pr-sequence of the inverse of a
nonsingular real symmetric matrix $B$.
The next result demonstrates that, when $n \geq 3$,
$r'_0$ can always be determined from $\pr(B)$ if this sequence does not end with $111$.

\begin{rem}\label{Inverse with 01 and 011}
\rm{Let $n \geq 3$,
$B$ be a nonsingular real symmetric matrix with
$\pr(B) = r_0]r_1 \cdots r_{n-1} 1$ and
$r'_0$ be the 0th term of $\pr(B^{-1})$.

  \begin{enumerate}[label=(\roman*)]

   \item If $r_{n-1}r_n = 01$, then
       $r'_0 = 1$.

   \item If $r_{n-2}r_{n-1}r_n = 011$, then
       $r'_0 = 0$.

  \end{enumerate}

(i) is immediate from Theorem \ref{Inverse pr Theorem},
since $B$ obviously has a
principal minor of order $n-1$ that is zero.
As for (ii), first, notice that the penultimate term of $\epr(B)$ must be $\tt{A}$, as $\tt{NSA}$ is prohibited; therefore, $B$ does not have a principal minor of order $n-1$ that is zero, implying that $r'_0 = 0$.
}
\end{rem}

The next proposition generalizes
a particular case of \cite[Lemma 4.5]{ORIGINAL}.

\begin{prop}
Let $B$ be a real symmetric matrix with
$\pr(B) = r_0]r_1 \cdots r_n$.
Suppose that $\pr(B)$ does not contain three consecutive $1$s and that $r_0]r_1 \neq 1]1$.
Then, for any $m$ with $1 \leq m \leq n$,
there exists a principal submatrix $B'$ of
$B$ such that $\pr(B') = r_0]r_1 \cdots r_m$.
\end{prop}

\bpf
Let $1 \leq m \leq n$.
By \cite[Lemma 4.5]{ORIGINAL}, $B$ has a principal submatrix
$B'$ with $\pr(B') = r'_0]r_1r_2 \cdots r_m$.
Since $B$ does not contain both a zero and
a nonzero diagonal entry, it follows that
$r'_0]r_1 = r_0]r_1$, and therefore
$\pr(B') = r_0]r_1 \cdots r_m$.
\epf

\begin{cor}\label{0r1r2...rn with no 111, etc.}
Let $\sigma = r_0]r_1 \cdots r_n$ be a pr-sequence
not containing three consecutive $1$s.
Suppose $r_0]r_1 \neq 1]1$.
If any initial subsequence of $\sigma$ is
unattainable, then $\sigma$ is unattainable.
\end{cor}

It was shown in \cite{ORIGINAL} that appending $0$ to
the end of an attainable pr-sequence results in
a new attainable pr-sequence;
but, what if 0 is appended to an
unattainable pr-sequence?
For example, if we append 0 to $1]1011$, which is
unattainable (see \cite[Table 5.4]{ORIGINAL}),
we obtain the attainable pr-sequence $1]10110$
(see \cite[Table 6.1]{ORIGINAL}).
However, there are some cases
where appending $0$ preserves unattainability.
The next observation, a consequence of
Corollary \ref{0r1r2...rn with no 111, etc.}, illustrates this.

\begin{obs}\label{Obs.-Append-zeros-to-unattain.-sequences}
Let $r_0]r_1 \cdots r_n$ be an unattainable pr-sequence
not containing three consecutive $1$s.
Suppose $r_0]r_1 \neq 1]1$.
Then $r_0]r_1 \cdots r_n 0$ is also unattainable.
\end{obs}

Propositions \ref{ANSNSN(SN)(N)} and
\ref{SNSNSN(SN)(N)} below are corollaries to
Theorem \ref{ANSNSN...A}.

\begin{prop}\label{ANSNSN(SN)(N)}
Let $B$ be a real symmetric matrix with
$\epr(B) = \tt{ANSNSN} \cdots$.
Then, for $k \geq 1$, $\ell_{2k} = \tt{N}$.
Furthermore, $\epr(B) = \tt{ANSNSN\overline{SN}}\hspace{0.04cm}\overline{\tt{N}}$
or
$\epr(B) = \tt{ANSNSN\overline{SN}A}$.
\end{prop}

\bpf
Let $k \geq 1$.
By hypothesis, the first assertion holds for $k \leq 3$.
Suppose $\ell_{2k} \neq \tt{N}$ for some $k >3$.
By the Inheritance Theorem, $B$ has a nonsingular
$2k \times 2k$ principal submatrix with epr-sequence
$\tt{ANXNYN} \cdots \tt{A}$, where
$\tt{X}$, $\tt{Y}$, $\tt{Z} \in
\{\tt{A}, \tt{S}, \tt{N}\}$.
By the $\tt{NN}$ Theorem,
$\tt{X}$ and $\tt{Y}$ are not $\tt{N}$.
Since $\tt{NAN}$ is prohibited,
$\tt{X} = \tt{Y} = \tt{S}$,
a contradiction to Theorem \ref{ANSNSN...A}.
%For the final statement, notice that $\ell_{2k+1}$ is
%preceded and followed by $\tt{N}$, and hence
%it cannot be $\tt{A}$, as $\tt{NAN}$ is prohibited;
%the assertion now follows from the $\tt{NN}$ Theorem.
The final assertion is immediate from
the $\tt NN$ Theorem and the fact that
$\tt NAN$ is prohibited.
\epf

\begin{cor}\label{ANSNSN(SN)(N) Corollary}
The pr-sequence
$0]1 01 01 01 \overline{01} 1 \overline{0}$
is not attainable by a real symmetric matrix.
\end{cor}

\bpf
Since $0]1 01 01 01 \overline{01} 1$
satisfies the hypothesis of
Observation \ref{Obs.-Append-zeros-to-unattain.-sequences},
it suffices to show that this sequence is not attainable.
% It was not necessary to use that observation.
Suppose that there is a real symmetric matrix $B$ with
$\pr(B) = 0]1 01 01 01 \overline{01} 1$ and
$\epr(B) = \ell_1 \ell_2 \cdots \ell_n$.
Obviously, $\ell_1 = \ell_n =\tt{A}$ and
$\ell_2 = \ell_4 = \ell_6= \tt{N}$.
Since $\tt{NAN}$ is prohibited,
$\ell_3 = \ell_5 = \tt{S}$.
Hence,
$\epr(B) = \tt{ANSNSN \cdots \tt{XA}}$,
where $\tt{X}$ is not $\tt{N}$,
which contradicts Proposition \ref{ANSNSN(SN)(N)}.
\epf

\begin{prop}\label{SNSNSN(SN)(N)}
Let $B$ be a real symmetric matrix with
$\epr(B) = \tt{SNSNSN} \cdots$.
Then, for $k \geq 1$, $\ell_{2k} = \tt{N}$.
Furthermore, $\epr(B) = \tt{SNSNSN\overline{SN}}\hspace{0.04cm}\overline{\tt{N}}$
or
$\epr(B) = \tt{SNSNSN\overline{SN}A}$.
\end{prop}

\bpf
Let $k \geq 1$.
By hypothesis, the first assertion holds for $k \leq 3$.
Suppose $\ell_{2k} \neq \tt{N}$ for some $k >3$.
By the Inheritance Theorem, $B$ has a nonsingular
$2k \times 2k$ principal submatrix with an epr-sequence
$\tt{XNYNZN} \cdots \tt{A}$, where
$\tt{X}$, $\tt{Y}$, $\tt{Z} \in
\{\tt{A}, \tt{S}, \tt{N}\}$.
By the $\tt{NN}$ Theorem,
$\tt{X}$, $\tt{Y}$ and $\tt{Z}$ are not $\tt{N}$.
Since $\tt{NAN}$ is prohibited, $\tt{Y} = \tt{Z} = \tt{S}$.
Since $\tt{SN} \cdots \tt{A} \cdots$ is prohibited,
$\tt{X} \neq \tt{S}$, and hence $\tt{X} = \tt{A}$,
a contradiction to Theorem \ref{ANSNSN...A}.
%%%For the final statement, notice that $\ell_{2k+1}$ is
%%%preceded and followed by an $\tt{N}$, and hence
%%%it cannot be $\tt{A}$, as $\tt{NAN}$ is prohibited;
%%%the assertion follows from the $\tt{NN}$ Theorem.
As in Proposition \ref{ANSNSN(SN)(N)},
the final assertion follows from
the $\tt{NN}$ Theorem and
the fact that $\tt{NAN}$ is prohibited.
\epf

\begin{cor}\label{SNSNSN(SN)(N) Corollary}
The pr-sequence
$1]1 01 01 01 \overline{01} 1 0 \overline{0}$
is not attainable by a real symmetric matrix.
\end{cor}

\bpf
Suppose there is a real symmetric matrix $B$ with
$\pr(B) = 1]1 01 01 01 \overline{01} 1 0 \overline{0}$.
Let $\epr(B) = \ell_1 \ell_2 \cdots \ell_n$.
Obviously, $\ell_1 = \tt{S}$ and
$\ell_2 = \ell_4 = \ell_6= \tt{N}$.
Since $\tt{NAN}$ is prohibited,
$\ell_3 = \ell_5 = \tt{S}$.
Hence,
$\epr(B) = \tt{SNSNSN \cdots \tt{XYN} \overline{\tt{N}}}$,
where $\tt{X}$ and $\tt{Y}$ are both not $\tt{N}$,
which contradicts Proposition \ref{SNSNSN(SN)(N)}.
\epf

Before proving the main result of this section,
we need a lemma.

%%%\begin{lem}
%%%Let $n \geq 4$ and
%%%$B$ be a real symmetric matrix with
%%%$\pr(B) = r_0]r_1 \cdots r_n$.
%%%Suppose $r_1r_2 \cdots r_n$ does not
%%%contain three consecutive $1$s.
%%%Then, for all $k$ such that $2 \leq k \leq \rank(B)-2$,
%%%$r_k r_{k+1}$ contains $0$.
%%%\end{lem}
%%%
%%%\bpf
%%%Let $2 \leq k \leq \rank(B)-2$.
%%%Suppose to the contrary that $r_k r_{k+1} = 11$.
%%%By hypothesis, $r_{k-1} = r_{k+2} = 0$;
%%%by the $0110$ Theorem, $r_j = 0$ for $j \geq k+2$,
%%%implying that $\rank(B) < k+2$, a contradiction.
%%%\epf

\begin{lem}\label{No 111 lemma}
Let $B$ be a real symmetric matrix
with $\pr(B) = r_0]r_1 \cdots r_n$.
Suppose $r_1r_2 \cdots r_n$
does not contain three consecutive $1$s.
Let $1 \leq k \leq \rank(B)-2$.
If $r_k r_{k+1} = 01$, then
either
$r_{k+2}r_{k+3}\cdots r_n = \overline{01} 1 \overline{0}$
or
$r_{k+2}r_{k+3}\cdots r_n =
\overline{01}01\overline{0}$
\end{lem}

\bpf
Suppose $r_k r_{k+1} = 01$.
We proceed by examining the only
two possibilities for $r_{k+2}$.

\textit{Case 1}: $r_{k+2} = 1$.
Now we have $r_k r_{k+1}r_{k+2} = 011$.
If $n=k+2$, then we are done.
Now, suppose $n>k+2$.
By hypothesis, $r_{k+3} = 0$, and therefore,
by the $0110$ Theorem,
$r_{k+2}r_{k+3} \cdots r_n = 1 \overline{0}$, where
$\overline{0}$ is non-empty.

\textit{Case 2}: $r_{k+2} = 0$.
Now we have $r_k r_{k+1}r_{k+2} = 010$.
Then, as $\rank(B) \geq k+2$,
by the $00$ Theorem, $r_{k+3} \neq 0$;
hence, $r_{k+3} = 1$, and so
$r_{k+2} r_{k+3} = 01$.
If $n=k+3$, then we are done.
Suppose $n > k+3$.
If $\rank(B) = k+3$, then we have
$r_{k+2}r_{k+3}\cdots r_n = 01\overline{0}$,
where $\overline{0}$ is non-empty.
Suppose $\rank(B) > k+3$, i.e.,
suppose $\rank(B) \geq k+4$.
Thus, so far we have
$r_k r_{k+1} r_{k+2} r_{k+3} = 0101$, where
$r_{k+2} r_{k+3} = 01$ and
$1 \leq k+2 \leq \rank(B)-2$.
Since $n$ is finite, it is evident that
reimplementing the steps above by replacing
$k$ with $k+2$, and repeating this process until reaching
the last term of the sequence, yields the desired conclusion.
\epf

With the next theorem, we classify all
the attainable pr-sequences of order $n \geq 3$ not
containing three consecutive 1s.

\begin{thm}\label{No-111 characterization}
Let $n \geq 3$.
A pr-sequence of order $n$ not
containing three consecutive $1$s is attainable by
a real symmetric matrix
if and only if
it is one of the following sequences.

\begin{description}

   \item[\textnormal{1.}]
        $0]100\overline{0}$.

   \item[\textnormal{2.}]
        $0]1 \overline{01}01\overline{0}$.

   \item[\textnormal{3.}]
        $0]1 01 1 \overline{0}$.

   \item[\textnormal{4.}]
        $0]1 01 01 1 \overline{0}$.

   \item[\textnormal{5.}]
        $0]1 1 0\overline{0}$.

   \item[\textnormal{6.}] $0]1101\overline{0}$.

   \item[\textnormal{7.}] $0]11011\overline{0}$.

   \item[\textnormal{8.}] $1]000\overline{0}$.

   \item[\textnormal{9.}] $1]01 0\overline{0}$.

   \item[\textnormal{10.}]
        $1] 01 \overline{01}01\overline{0}$.

   \item[\textnormal{11.}] $1]01\overline{01}1\overline{0}$.

   \item[\textnormal{12.}] $1]100\overline{0}$.

   \item[\textnormal{13.}]
       $1]1 \overline{01}01 0\overline{0}$.

   \item[\textnormal{14.}] $1]1 01 1 0\overline{0}$.

   \item[\textnormal{15.}] $1]1 01 01 1 0\overline{0}$.

  \end{description}
\end{thm}

\bpf
Let $B$ be a real symmetric matrix with
$\pr(B) = r_0]r_1 \cdots r_n$ not
containing three consecutive $1$s.
Since $0]0 \cdots$ is forbidden by definition,
$r_0]r_1 \in \{0]1, 1]0, 1]1 \}$.
We proceed by examining all the
possibilities for $r_0]r_1r_2$.

\textit{Case i}: $r_0]r_1r_2 = 0]10$. If $r_3 = 0$, then,
by the $00$ Theorem, we have sequence (1).
Suppose $r_3 = 1$.
Hence, $\pr(B)$ starts $0]101 \cdots$.
If $\rank(B) = 3$, then
$\pr(B) = 0]101 \overline{0}$,
which is sequence (2).
Now, suppose $\rank(B) > 3$.
Then $r_2 r_3 = 01$ and $1 \leq 2 \leq \rank(B) - 2$;
hence, by applying Lemma \ref{No 111 lemma} to $\pr(B)$,
starting with $k=2$, we have either
$\pr(B)=0]1 01 \overline{01} 01 \overline{0}$ or
$\pr(B)=0]1 01 \overline{01} 1 \overline{0}$.
Hence, by Corollary \ref{ANSNSN(SN)(N) Corollary},
$\pr(B)$ is one of the sequences (2), (3) and (4).

\textit{Case ii}:  $r_0]r_1r_2 = 0]11$.
By hypothesis, $r_3 = 0$.
If $\rank(B) = 2$, then
$\pr(B) = 0]110 \overline{0}$,
which is sequence (5).
Now suppose $\rank(B) > 2$.
Then $n >3$ and, by the $00$ Theorem, $r_4 \neq 0$,
implying that $r_4 = 1$.
Hence, $\pr(B)$ starts $0]1101 \cdots$.
If $n=4$, then we have sequence (6).
Suppose $n>4$.
If $r_5 = 1$, then, by the $0110$ Theorem,
we must have sequence (7), where
$\overline{0}$ may be empty.
Now, suppose $r_5 = 0$.
If $n = 5$, then we have sequence (6).
Suppose $n>5$.
Thus far we have $\pr(B) = 0]11010 \cdots$;
it follows from \cite[Theorem 7.2]{ORIGINAL} that
$r_6 = 0$, and therefore, by the $00$ Theorem,
we have sequence (6).

\textit{Case iii}:  $r_0]r_1 = 1]0$.
If $r_2 = 0$, then, by the $00$ Theorem,
we have sequence (8).
Now, suppose $r_2 = 1$.
Hence, $\pr(B)$ starts $1]01 \cdots$.
If $\rank(B) = 2$, then
$\pr(B) = 1]01 0\overline{0}$,
which is sequence (9).
Now, suppose $\rank(B) >2$.
Then $r_1r_2 = 01$ and $1 \leq 1 \leq \rank (B) - 2$;
hence, by applying Lemma \ref{No 111 lemma} to $\pr(B)$,
starting with $k=1$,
we have either
$\pr(B)=1]01 \overline{01} 01 \overline{0}$ or
$\pr(B)=1]01 \overline{01}1\overline{0}$.
Thus, $\pr(B)$ is either sequence (10) or (11).

\textit{Case iv}:  $r_0]r_1 = 1]1$.
By hypothesis, $r_2 = 0$.
If $r_3 = 0$, then the $00$ Theorem implies that
we have sequence (12).
Now, suppose $r_3 = 1$.
Hence, $\pr(B)$ starts $1]101 \cdots$.
Suppose $\rank(B) = 3$;
then $\pr(B) = 1]1 01\overline{0}$, and,
by \cite[Theorem 4.1]{ORIGINAL},
$\overline{0}$ is non-empty, implying that
$\pr(B) = 1]1 01 0\overline{0}$, which is
sequence (13).
Now, suppose $\rank(B) >3$.
Then $r_2r_3 = 01$ and $1 \leq 2 \leq \rank (B) - 2$;
hence, by applying Lemma \ref{No 111 lemma} to $\pr(B)$,
starting with $k=2$, we have
$\pr(B)=1]1 01 \overline{01} 01 \overline{0}$ or
$\pr(B)=1]1 01 \overline{01} 1 \overline{0}$;
again, it follows from \cite[Theorem 4.1]{ORIGINAL} that
in either case $\overline{0}$ must be non-empty,
and therefore
$\pr(B) = 1]1 01 \overline{01} 01 0\overline{0}$ or
$\pr(B) = 1]1 01 \overline{01} 1 0\overline{0}$.
Hence, by
Corollary \ref{SNSNSN(SN)(N) Corollary},
$\pr(B)$ is one of the sequences (13), (14) and (15).

For the other direction, since appending $0$ to
the end of an attainable sequence results in
another attainable sequence
(see \cite[Theorem 2.6]{ORIGINAL}),
it suffices to establish the attainability of
each sequence when $\overline{0}$ is empty.
We assume that the sequence under
consideration has order $n \geq 3$ and provide an
$n\times n$ real symmetric matrix that attains it.
\begin{description}
\item[1.]$0]100\overline{0}$:
$\pr(J_3) = 0]100$.

\item[2.]$0]1 \overline{01}01\overline{0}$:
$\pr((A(C_n))^{-1}) = 0]1 \overline{01}01$,
with $n$ odd (see \cite[Lemma 3.4]{ORIGINAL} and
Remark \ref{Inverse with 01 and 011}).
%% also, see \cite{FIELDS}, Section \ref{s:restric}.3, item 2.

\item[3.]$0]1 01 1\overline{0}$:
$\pr(J_4-2I_4) = 0]1 01 1$.

\item[4.]$0]1 01 01 1\overline{0}$:
$\pr(M_{0101011}) = 0]1 01 01 1$, where
$M_{0101011}$ appears in \cite[p. 2153]{ORIGINAL}.

\item[5.]$0]1 1 0 \overline{0}$:
$\pr(J_1 \oplus J_2) = 0]1 1 0$.

\item[6.] $0]1101 \overline{0}$:
$\pr(J_4-3I_4) = 0]1101$.

\item[7.] $0]11011 \overline{0}$:
$\pr(J_5-3I_5) = 0]11011$.

\item[8.] $1]000\overline{0}$:
$\pr(0_3) = 1]000$.

\item[9.] $1]01 0\overline{0}$:
$\pr((J_2 - I_2) \oplus 0_1) = 1]01 0$.

\item[10.]$1]01 \overline{01}01\overline{0}$:
$\pr(A(P_n)) = 1] 01 \overline{01}01$, with $n$ even
(see \cite[Lemma 3.3]{ORIGINAL}).

\item[11.] $1]01\overline{01}1\overline{0}$:
$\pr(A(C_n)) = 1]01\overline{01}1$, with $n$ odd
(see \cite[Lemma 3.4]{ORIGINAL}).

\item[12.] $1]100\overline{0}$:
$\pr(J_1 \oplus 0_2) = 1]100$.

\item[13.] $1]1 \overline{01}01 0\overline{0}$:
$\pr((A(C_{n-1}))^{-1} \oplus 0_1) =
1]1\overline{01} 01 0$, with $n$ even
(see \cite[Lemma 3.4]{ORIGINAL},
     Remark \ref{Inverse with 01 and 011} and
     \cite[Theorem 2.3]{ORIGINAL}).
%%% This also follows in the following manner: applying
%%% Remark \ref{Inverse with 01 and 011} to the attainable
%%% sequence $1]0\overline{10}11$ (of appropriate order) and
%%% later applying \cite[Theorem 2.3]{ORIGINAL} will result
%%% in the desired sequence.

\item[14.] $1]1 01 1 0\overline{0}$:
$\pr((J_4 - 2I_4) \oplus 0_1) = 1]1 01 1 0$.

\item[15.] $1]1 01 01 1 0\overline{0}$:
$\pr( M_{0101011} \oplus 0_1) = 1]1 01 01 1 0$,
where $M_{0101011}$
appears in \cite[p. 2153]{ORIGINAL}.
%%%(apply  \cite[Theorem 2.3]{ORIGINAL} to the matrix
%%%$M_{0101011}$ in \cite[p. 2153]{ORIGINAL}).

\end{description}
That concludes the proof.
\epf

We conclude this section with a classification of
the attainable pr-sequences that
only contain three consecutive $1$s in the initial
subsequence 1]11.
The primary motivation for including this result is
its application in Section \ref{s:N in length 3}.

\begin{prop}\label{SS(NS)SN(N) and SS(NS)NAA}
The epr-sequences
$\tt{SS NS NS\overline{NS}SN\overline{N}}$ and
$\tt{SS NS\overline{NS}NAA}$
are not attainable by a real symmetric matrix.
\end{prop}

\bpf
Suppose to the contrary that there is a
real symmetric matrix $B$ with
$\epr(B) = \tt{SSNSNS\overline{NS}SN\overline{N}}$.
Notice that $\rank(B)$ is odd.
If $\tt{\overline{NS}}$ is empty, then
we have a contradiction to Proposition \ref{SSNSNSS}.
So, suppose $\tt{\overline{NS}}$ is non-empty.
Let $B[\alpha]$ be a nonsingular $1 \times 1$
principal submatrix of $B$.
By the Schur Complement Theorem,
$\rank(B/B[\alpha])$ is even,
$\rank(B/B[\alpha]) \geq 8$,
and $\epr(B/B[\alpha]) = \tt{XNYNZN} \cdots$,
where $\tt{X}$, $\tt{Y}$, $\tt{Z} \in
\{\tt{A}, \tt{S}, \tt{N}\}$.
Then, as $\rank(B/B[\alpha]) \geq 8$,
by the $\tt{NN}$ Theorem,
$\tt{X}$, $\tt{Y}$ and $\tt{Z}$ are not $\tt{N}$.
Since $\tt{NAN}$ is prohibited, $\tt{Y} = \tt{Z} = \tt{S}$.
Thus, we have $\epr(B/B[\alpha]) = \tt{XNSNSN} \cdots$,
where $\tt{X}$ is not $\tt{N}$.
It follows from Propositions \ref{ANSNSN(SN)(N)} and
\ref{SNSNSN(SN)(N)} that $\rank(B/B[\alpha])$ is odd,
a contradiction.

Now, suppose
$\tt{SSNS\overline{NS}NAA}$ is attainable.
Then applying \cite[Observation 2.19(2)]{EPR} to this
sequence implies that
$\tt{SSNSNS\overline{NS}SN}$ is attainable,
a contradiction to the first assertion.
\epf

\begin{cor}\label{SS(NS)SN(N) and SS(NS)NAA Corollary}
The pr-sequence
$1]11 01 01  \overline{01} 1 \overline{0}$
is not attainable by a real symmetric matrix.
\end{cor}

\bpf
Suppose that there is a real symmetric matrix $B$ with
$\pr(B) = 1]11 01 01  \overline{01} 1 \overline{0}$ and
$\epr(B) = \ell_1 \ell_2 \cdots \ell_n$.
Obviously, $\ell_1 = \tt{S}$ and
$\ell_3 = \ell_5 = \tt{N}$.
By the $\tt{NN}$ Theorem,
and because $\tt{NAN}$ is prohibited,
$\ell_4 = \tt{S}$.
Since $\ell_2$ is not $\tt{N}$,
it follows from Proposition \ref{SAN} that
$\ell_2 = \tt{S}$.
Hence,
$\epr(B) = \tt{SS NS N \cdots}$.
We examine two cases.

\textit{Case 1}: $\overline{0}$ is empty.
Notice that
$\pr(B) = 1]11 01 01  \overline{01} 1 =
          1]11 01  \overline{01} 01 1$.
Moreover, $\ell_n = \tt{A}$ and
$\ell_{i} = \tt{N}$ for all odd $i$ with
$3 \leq i \leq n-2$.
Then, as $\tt{NAN}$ is prohibited,
$\ell_{j} = \tt{S}$ for all even $j$ with
$4 \leq j \leq n-3$.
Therefore, we have
$\epr(B) = \tt{SS NS \overline{NS} NXA}$,
where $\tt{X}$ is not $\tt{N}$.
Since $\tt{NSA}$ is prohibited,
$\tt{X} = \tt{A}$, which contradicts
Proposition \ref{SS(NS)SN(N) and SS(NS)NAA}.

\textit{Case 2}: $\overline{0}$ is non-empty.
Thus,
$\pr(B) = 1]11 01 01  \overline{01} 1 0\overline{0} =
          1]11 01  \overline{01} 01 1 0\overline{0}$.
As in the preceding case, the fact that
$\tt{NAN}$ is prohibited implies that
$\epr(B) =
\tt{SS NS \overline{NS} NX Y N\overline{N}}$,
where $\tt{X}$ and $\tt{Y}$ are not $\tt{N}$.
By Theorem \ref{SNA}, $\tt{X} = \tt{S}$.
Then, as $\tt{NSA}$ is prohibited, $\tt{Y} = \tt{S}$.
Hence,
$\epr(B) = \tt{SS NS NS\overline{NS}SN\overline{N}}$,
a contradiction to
Proposition \ref{SS(NS)SN(N) and SS(NS)NAA}.
\epf

\begin{prop}\label{No-111 characterization 2}
Let $n \geq 3$.
A pr-sequence $r_0]r_1 \cdots r_n$, with
$r_1 r_2 \cdots r_n$ not containing three consecutive $1$s,
is attainable by a real symmetric matrix
if and only if
it is one of the sequences in
Theorem \ref{No-111 characterization} or
one of the following sequences.

\begin{description}

   %16
   \item[\textnormal{16.}]
        $1]110\overline{0}$.

   %17
   \item[\textnormal{17.}]
        $1]11 \overline{01} 01 \overline{0}$.

   %18
   \item[\textnormal{18.}]
        $1]11 01 1 \overline{0}$.

\end{description}

\end{prop}

\bpf
Let $B$ be a real symmetric matrix with
$\pr(B) = r_0]r_1 \cdots r_n$.
Suppose $r_1 r_2 \cdots r_n$ does not
contain three consecutive $1$s.
If $r_0]r_1r_2 \neq 1]11$, then $\pr(B)$ does not contain three consecutive $1$s, and therefore it is one of the sequences listed in Theorem \ref{No-111 characterization}.
Thus, suppose $r_0]r_1r_2 = 1]11$.
By hypothesis, $r_3 = 0$.
If $n=3$, then $\pr(B)$ is sequence (16).
So, suppose $n>3$.
If $r_4 = 0$, then, by the $00$ Theorem,
$\pr(B)$ is sequence (16).
Now, suppose $r_4 = 1$.
Then $\pr(B)$ starts $1]1101 \cdots$.
If $\rank(B) = 4$, then
$\pr(B) = 1]1101\overline{0}$,
which is sequence (17).
Now, suppose $\rank(B) > 4$.
Hence, $r_3r_4 = 01$ and $1 \leq 3 \leq \rank(B) - 2$.
It follows from applying Lemma \ref{No 111 lemma} to
$\pr(B)$, starting with $k = 3$, that
$\pr(B)=1]11 01 \overline{01} 01\overline{0}$ or
$\pr(B)=1]11 01\overline{01} 1 \overline{0}$.
Hence, by
Corollary \ref{SS(NS)SN(N) and SS(NS)NAA Corollary},
$\pr(B)$ is either sequence (17) or sequence (18).\\
\indent
For the other direction, as in
Theorem \ref{No-111 characterization},
it suffices to show that each sequence is attainable when
$\overline{0}$ is empty.
By \cite[Theorem 3.7]{ORIGINAL}, the sequences
$1]110$ and $1]11 01 1$ are attainable
by $Q_{3,1}$ and $Q_{5,1}$, respectively.
Finally, $1]11 \overline{01} 01$ is attained by
$(A(F_n))^{-1}$ (see \cite[Lemma 3.5]{ORIGINAL}),
where $n$ is even and $F_n$ is
the graph on $n$ vertices formed by adding a
pendent edge to $C_{n-1}$.
\epf

\section{Epr-sequences with an $\tt{N}$ in every subsequence of length 3}\label{s:N in length 3}
$\indent$
This section focuses on epr-sequences with an $\tt{N}$ in every subsequence of length $3$, and culminates with
a classification of all the attainable epr-sequences
with this property.

The sequence accounted for in the next result is of particular relevance to the main result at the end of this section.

\begin{prop}\label{NAAN}
Let $n \geq 3$ and
$B=[b_{ij}]$ be the $n \times n$ real symmetric matrix with
$b_{ij} = (i-j)^2$.
Then $\epr(B) = \tt{NAA\overline{N}}$.
\end{prop}

\bpf
Suppose that $\epr(B) = \ell_{1}\ell_{2}\cdots \ell_{n}$.
It is easy to verify the assertion for $n=3$.
Suppose $n > 3$.
Obviously, $\ell_{1} = \tt{N}$.
Let $p,q,r,s \in \{1,2, \dots, n\}$,
where $p < q < r < s$.
Since every off-diagonal entry of $B$ is nonzero, we have
$B_{pq} = -(b_{pq})^2 \neq 0$ and
$B_{pqr} = 2b_{pq}b_{pr}b_{qr} \neq 0$.
A simple computation reveals that
the order-$4$ principal minor $B_{pqrs}$ is given by
\[(b_{ps}b_{qr})^2 + (b_{pr}b_{qs})^2 + (b_{pq}b_{rs})^2 -
2b_{pr}b_{ps}b_{qr}b_{qs} - 2b_{pq}b_{ps}b_{qr}b_{rs} -
2b_{pq}b_{pr}b_{qs}b_{rs}= \]
\[((p-s)(q-r))^4 + ((p-r)(q-s))^4 + ((p-q)(r-s))^4\]
\[- 2((p-r)(p-s)(q-r)(q-s))^2 - 2((p-q)(p-s)(q-r)(r-s))^2\]
\[- 2((p-q)(p-r)(q-s)(r-s))^2 = 0.\]
Hence, we have $\ell_{2} = \ell_{3} = \tt{A}$ and
$\ell_{4} = \tt{N}$. The conclusion now follows from
Proposition \ref{NXXN}.
\epf

\begin{obs}\label{No 111 implies an N in every subseq. of lenght 3}
If an attainable pr-sequence
does not contain three consecutive 1s,
then an attainable epr-sequence associated with
it contains an $\tt{N}$ in every subsequence of length $3$.
\end{obs}

\begin{rem}\label{Converse of No-111 implies N in every order 3-subs}
\rm{The converse of Observation
\ref{No 111 implies an N in every subseq. of lenght 3} is false.
An attainable epr-sequence starting $\tt{SS} \cdots$, or
starting $\tt{SA} \cdots$, with an $\tt{N}$ in every
subsequence of length $3$, provides a counterexample.
It can be deduced that
all counterexamples are of that form,
and therefore that the converse of Observation
\ref{No 111 implies an N in every subseq. of lenght 3}
is true if additionally we assume that
the pr-sequence does not start with $1]11$.
}
\end{rem}

\begin{obs}\label{111 only at the start}
Let $n \geq 3$ and $B$ be a real symmetric matrix with
$\pr(B) = r_0]r_1 \cdots r_n$.
Then
$\epr(B)$ contains an $\tt{N}$ in
every subsequence of length $3$
if and only if
$r_1 r_2 \cdots r_n$ does not contain three consecutive $1$s
\end{obs}

Observation \ref{111 only at the start}
suggests that we can use
Theorem \ref{No-111 characterization} and
Proposition \ref{No-111 characterization 2} to
classify all the epr-sequences with an $\tt{N}$ in
every subsequence of length $3$,
as the pr-sequences associated with these epr-sequences
must be those listed on these results.

\begin{thm}\label{N in every 3 characterization}
Let $n \geq 3$.
An epr-sequence of order $n$
with an $\tt{N}$ in every subsequence of length $3$ is
attainable by a real symmetric matrix
if and only if
it is one of the following sequences.

\begin{description}

   \item[\textnormal{1.}]
       $\tt{ANN\overline{N}}$.

   \item[\textnormal{2a.}]
       $\tt{A\overline{NS}NA}$.

   \item[\textnormal{2b.}]
       $\tt{A \overline{NS}NS N\overline{N}}$.

   \item[\textnormal{3a.}]
       $\tt{ANAA}$.

   \item[\textnormal{3b.}]
       $\tt{ANSSN\overline{N}}$.

   \item[\textnormal{4a.}]
       $\tt{ANSNAA}$.

   \item[\textnormal{4b.}]
       $\tt{A NS NSSN\overline{N}}$.

   \item[\textnormal{5a.}]
       $\tt{AAN\overline{N}}$.

   \item[\textnormal{5b.}]
       $\tt{ASN\overline{N}}$.

   \item[\textnormal{6a.}]
       $\tt{AANA}$.

   \item[\textnormal{6b.}]
       $\tt{ASNSN\overline{N}}$.

   \item[\textnormal{7a.}]
       $\tt{AANAA}$.

   \item[\textnormal{7b.}]
       $\tt{ASNSSN\overline{N}}$.

   \item[\textnormal{8.}]
       $\tt{NNN\overline{N}}$.

   \item[\textnormal{9.}]
       $\tt{NSN\overline{N}}$.

   \item[\textnormal{10a.}]
       $\tt{NS\overline{NS}NA}$.

   \item[\textnormal{10b.}]
       $\tt{NS\overline{NS}NS N\overline{N}}$.

   \item[\textnormal{11a.}]
       $\tt{N\overline{SN}AA}$.

   \item[\textnormal{11b.}]
       $\tt{N \overline{SN} SSN\overline{N}}$.

   \item[\textnormal{11c.}]
       $\tt{NAAN\overline{N}}$.

   \item[\textnormal{12.}]
       $\tt{SNN\overline{N}}$.

   \item[\textnormal{13.}]
       $\tt{S\overline{NS}NS N\overline{N}}$.

   \item[\textnormal{14.}]
       $\tt{SNSSN\overline{N}}$.

   \item[\textnormal{15.}]
       $\tt{S NS NSSN\overline{N}}$.

   \item[\textnormal{16a.}]
       $\tt{SAN\overline{N}}$.

   \item[\textnormal{16b.}]
       $\tt{SSN\overline{N}}$.

   \item[\textnormal{17a.}]
       $\tt{SS\overline{NS}NA}$.

   \item[\textnormal{17b.}]
       $\tt{SS\overline{NS}NS N\overline{N}}$.

   \item[\textnormal{18a.}]
       $\tt{SSNAA}$.

   \item[\textnormal{18b.}]
       $\tt{SSNSSN\overline{N}}$.

  \end{description}
\end{thm}

\bpf
Let $B$ be a real symmetric matrix with
$\epr(B) = \ell_1 \ell_2 \cdots \ell_n$.
Suppose that $\epr(B)$ contains an $\tt{N}$ in every
subsequence of length $3$.
It follows from Observation \ref{111 only at the start}
that $\pr(B)$ is one of the sequences listed in
Theorem \ref{No-111 characterization} or
Proposition \ref{No-111 characterization 2}.
We examine the 18 possible cases.

\textit{Case 1}: $\pr(B) = 0]100\overline{0}$.
Obviously, $\epr(B) = \tt{ANN \overline{N}}$,
which is sequence (1).

\textit{Case 2}: $\pr(B) = 0]1 \overline{01}01\overline{0}$.
First, suppose $\overline{0}$ is empty.
Then, as $\tt{NAN}$ is prohibited,
$\epr(B) = \tt{A\overline{NS}NA}$,
which is sequence (2a).
Now, suppose $\overline{0}$ is non-empty.
Similarly, since $\tt{NAN}$ is prohibited,
$\epr(B) = \tt{A \overline{NS}NS N\overline{N}}$,
which is sequence (2b).

\textit{Case 3}: $\pr(B) = 0]1 01 1 \overline{0}$.
If $\overline{0}$ is empty, then,
as $\tt{NSA}$ is prohibited,
$\epr(B) = \tt{ANAA}$, which is sequence (3a).
If $\overline{0}$ is non-empty, then,
since $\tt{NSA}$ and $\tt{NAS}$ are prohibited,
we must have
$\tt{ANSSN\overline{N}}$ or
$\tt{ANAAN\overline{N}}$;
as the latter sequence is forbidden by
Theorem \ref{(A)ANAA(A)},
$\epr(B)$ is sequence (3b).

\textit{Case 4}: $\pr(B) = 0]1 01 01 1 \overline{0}$.
Suppose $\overline{0}$ is empty.
Since $\tt{NAN}$ and $\tt{NSA}$ are prohibited,
$\epr(B) = \tt{ANSNAA}$, which is sequence (4a).
Now suppose $\overline{0}$ is non-empty.
Then, as $\tt{NAN}$, $\tt{NAS}$ and $\tt{NSA}$ are
prohibited, $\epr(B)$ is either
$\tt{A NS NSSN\overline{N}}$ or
$\tt{A NS NAAN\overline{N}}$;
by Theorem \ref{SNA}, the latter sequence is forbidden,
and thus we have sequence (4b).

\textit{Case 5}: $\pr(B) = 0]1 1 0\overline{0}$.
Clearly,
$\epr(B) = \tt{AAN\overline{N}}$ or
$\epr(B) = \tt{ASN\overline{N}}$, which are
sequences (5a) and (5b), respectively.

\textit{Case 6}: $\pr(B) = 0]1101\overline{0}$.
If $\overline{0}$ is empty,
then, as $\tt{ASNA}$ is forbidden,
$\epr(B) = \tt{AANA}$, which is sequence (6a).
Suppose $\overline{0}$ is non-empty.
Since $\tt{NAN}$ is prohibited,
and because $\tt{ANS}$ must be initial,
$\epr(B) = \tt{ASNSN\overline{N}}$,
which is sequence (6b).

\textit{Case 7}: $\pr(B) = 0]11011\overline{0}$.
Suppose $\overline{0}$ is empty.
Since $\tt{NSA}$ and $\tt{ASN} \cdots \tt{A}$
are prohibited,
$\epr(B) = \tt{AANAA}$, which is sequence (7a).
Suppose $\overline{0}$ is non-empty.
Moreover, suppose $\ell_{2} = \tt{A}$.
Obviously, $\ell_n = \tt{N}$;
but, as $\tt{ANS}$ must be initial,
$\ell_4 = \tt{A}$, and therefore
Theorem \ref{(A)ANAA(A)} implies that
$\ell_n = \tt{A}$, a contradiction.
It follows that we must have $\ell_2 = \tt{S}$.
Since $\tt{ASN} \cdots \tt{A} \cdots$ is prohibited,
$\epr(B) = \tt{ASNSSN\overline{N}}$,
which is sequence (7b).

\textit{Case 8}: $\pr(B) = 1]000\overline{0}$.
Clearly, $\epr(B) = \tt{NNN\overline{N}}$,
which is sequence (8).

\textit{Case 9}: $\pr(B) = 1]01 0\overline{0}$.
Since $\tt{NAN}$ is prohibited,
$\epr(B) = \tt{NSN\overline{N}}$,
which is sequence (9).

\textit{Case 10}: $\pr(B) = 1] 01 \overline{01}01\overline{0}$.
If $\overline{0}$ is empty,
then, as $\tt{NAN}$ is prohibited,
$\epr(B) = \tt{NS\overline{NS}NA}$,
which is sequence (10a).
Similarly, if $\overline{0}$ is non-empty,
$\epr(B) = \tt{NS\overline{NS} NS N\overline{N}}$,
which is sequence (10b).

\textit{Case 11}: $\pr(B) = 1]01\overline{01}1\overline{0}$.
First, observe that
$\pr(B) = 1]0\overline{10}11\overline{0}$.
Suppose $\overline{0}$ is empty.
Since $\tt{NSA}$ and $\tt{NAN}$ are prohibited,
$\epr(B) = \tt{N\overline{SN}AA}$,
which is sequence (11a).
Suppose $\overline{0}$ is non-empty.
Moreover, suppose $\overline{10}$ is empty.
Then, as $\tt{NAS}$ and $\tt{NSA}$ are prohibited,
$\epr(B)$ is
$\tt{NSSN\overline{N}}$ or
$\tt{NAAN\overline{N}}$,
which are sequences (11b) and (11c), respectively.
Finally, suppose $\overline{10}$ is non-empty.
Since $\tt{NAS}$, $\tt{NSA}$ and $\tt{NAN}$ are prohibited,
$\epr(B)$ is either
$\tt{N SN\overline{SN} SSN\overline{N}}$ or
$\tt{N SN\overline{SN} AAN\overline{N}}$;
by Theorem \ref{SNA}, the latter sequence is forbidden,
and therefore $\epr(B)$ is sequence (11b),
with $\tt{\overline{SN}}$ non-empty.

\textit{Case 12}: $\pr(B) = 1]100\overline{0}$.
Obviously, $\epr(B) = \tt{SNN\overline{N}}$,
which is sequence (12).

\textit{Case 13}: $\pr(B) = 1]1 \overline{01}01 0\overline{0}$.
Since $\tt{SN} \cdots \tt{A} \cdots$ is prohibited,
it is immediate that
$\epr(B)=\tt{S\overline{NS}NS N\overline{N}}$,
which is sequence (13).

\textit{Case 14}: $\pr(B) = 1]1 01 1 0\overline{0}$.
As in Case 13, since
$\tt{SN} \cdots \tt{A} \cdots$ is prohibited, we must have
$\epr(B) = \tt{SNSSN\overline{N}}$,
which is sequence (14).

\textit{Case 15}: $\pr(B) = 1]1 01 01 1 0\overline{0}$.
Again, as
$\tt{SN} \cdots \tt{A} \cdots$ is prohibited,
we must have
$\epr(B) = \tt{S NS NSSN\overline{N}}$,
which is sequence (15).

\textit{Case 16}: $\pr(B) = 1]110\overline{0}$.
Clearly, $\epr(B)$ is either
$\tt{SAN\overline{N}}$ or $\tt{SSN\overline{N}}$,
which are sequences (16a) and (16b), respectively.

\textit{Case 17}: $\pr(B) = 1]11 \overline{01} 01 \overline{0}$.
Since $\tt{SAN} \cdots \tt{A} \cdots $ and
$\tt{SAN} \cdots \tt{S} \cdots$
are prohibited by Proposition \ref{SAN},
$\ell_{2} = \tt{S}$.
Suppose $\overline{0}$ is empty.
Then, as $\tt{NAN}$ is prohibited,
$\epr(B)=\tt{SS\overline{NS}NA}$, which is sequence (17a).
Suppose $\overline{0}$ is non-empty.
Similarly, since $\tt{NAN}$ is prohibited,
$\epr(B) = \tt{SS\overline{NS}NS N\overline{N}}$,
which is sequence (17b).

\textit{Case 18}: $\pr(B) = 1]11 01 1 \overline{0}$.
As in the preceding case, we must have
$\ell_{2} = \tt{S}$.
Suppose $\overline{0}$ is empty.
Since $\tt{NSA}$ is prohibited,
$\epr(B) = \tt{SSNAA}$,
which is sequence (18a).
Suppose $\overline{0}$ is non-empty.
Hence, the fact that
$\tt{NAS}$ and $\tt{NSA}$ are prohibited implies that
$\epr(B)$ is either
$\tt{SSNSSN\overline{N}}$ or
$\tt{SSNAAN\overline{N}}$;
by Theorem \ref{SNA}, the latter sequence is forbidden,
and thus $\epr(B)$ is sequence (18b).

For the other direction, we show that all
the sequences listed are attainable, and assume that
the sequence under consideration has order $n \geq 3$.
Sequence (1) is attained by $J_n$.
Sequence (2a) is attained by
$A((C_n)^{-1})$
(see \cite[Observation 3.1]{EPR} and the Inverse Theorem),
when $\tt{\overline{NS}}$ is non-empty, and by
\cite[Proposition 2.17]{EPR},
when $\tt{\overline{NS}}$ is empty.
As for (2b),
applying \cite[Observation 2.19(1)]{EPR} to (2a),
results in this sequence.
Sequence (3a) is attainable by
\cite[Proposition 2.17]{EPR}.
Sequence (3b) is attainable by applying
\cite[Observation 2.19(1)]{EPR} to (3a).
Sequence (4a) is attainable by \cite[Table 1]{EPR}, and
(4b) results from applying
\cite[Observation 2.19(1)]{EPR} to (4a).
Sequences (5a) and (5b) are attainable by
\cite[Theorem 4.6]{EPR}.
Sequence (6a) is attainable by \cite[Proposition 2.17]{EPR}, and
(6b) results from applying
\cite[Observation 2.19(1)]{EPR} to (6a).
Sequence (7a) is attainable by
\cite[Proposition 2.17]{EPR}, and
(7b) results from applying
\cite[Observation 2.19(1)]{EPR} to (7a).
Sequence (8) is attained by $0_n$.
As for (9), applying \cite[Observation 2.19(1)]{EPR}
to the sequence $\tt{NA}$, which is attained by
$J_2 - I_2$, results in this sequence.
Sequence (10a) is attainable by
\cite[Observation 3.1]{EPR}, and
(10b) results from applying
\cite[Observation 2.19(1)]{EPR} to (10a).
Sequence (11a) is attainable by
\cite[Observation 3.1]{EPR}, while
(11b) is obtained from applying
\cite[Observation 2.19(1)]{EPR} to (11a).
Sequence (11c) is attainable by Proposition \ref{NAAN}.
Sequence (12) is attainable by \cite[Theorem 4.6]{EPR}.
Sequences (13), (14) and (15) result from applying
\cite[Observation 2.19(2)]{EPR} to
(2a), (3a) and (4a), respectively.
Sequences (16a) and (16b) are attainable by
\cite[Theorem 4.6]{EPR}.
According to Proposition \ref{No-111 characterization 2},
the sequence $1]11 \overline{01} 01$ is attainable;
by Proposition \ref{SAN}, and because
$\tt{NAN}$ is prohibited, an attainable epr-sequence
associated with this pr-sequence, must be
$\tt{SS \overline{NS} NA}$, which is sequence (17a).
Sequence (17b) results from applying
\cite[Observation 2.19(2)]{EPR} to (17a).
Sequence (18a) is attainable by \cite[Table 5]{EPR}, and
(18b) is attainable by \cite[Corollary 2.20(2)]{EPR}.
\epf

If an epr-sequence is attainable, then the
pr-sequence associated with it must be attainable.
The converse is not true; this is because an epr-sequence
associated with a pr-sequence may not be unique,
since a 1 in the pr-sequence can correspond to an $\tt A$ or $\tt S$ in the epr-sequence.
For example, the epr-sequences $\tt NSSN$ and $\tt NAAN$,
which are associated with the pr-sequence $1]0110$,
are each attainable by a real symmetric matrix
(see \cite[Table 4]{EPR}).
We now show that, for real symmetric matrices, almost all attainable pr-sequences not containing three consecutive 1s are associated with a unique epr-sequence.

\begin{prop}\label{Uniqueness: general}
Let $n \geq 3$ and $\sigma$ be a pr-sequence that
is attainable by an $n \times n$ real symmetric matrix.
Suppose $\sigma$ does not contain three consecutive $1$s,
$\sigma  \neq 0]110  \OL{0}$ and that
$\sigma  \neq 1]0110 \OL{0}$.
Then there is a unique attainable epr-sequence associated with $\sigma$.
\end{prop}

\bpf
Since the attainable epr-sequences associated with pr-sequences not containing three consecutive 1s are
the epr-sequences (1a)--(15) listed in
Theorem \ref{N in every 3 characterization},
an attainable epr-sequence associated with
$\sigma$ must be one of these sequences.
Note that $\sigma$ is not associated with any of the epr-sequences (16a)--(18b), as these are the epr-sequences that are associated with
the pr-sequences listed in
Proposition \ref{No-111 characterization 2}.
We consider two cases.

\textit{Case 1}: $\sigma = 1]0 10 \OL{10} 11 0 \OL{0}$.
Observe that $\sigma$ is associated with
the epr-sequence (11b) in
Theorem \ref{N in every 3 characterization},
with $\tt \OL{SN}$ non-empty.
It is easy to see that
$\sigma$ is not associated with
any of the other epr-sequences listed in
Theorem \ref{N in every 3 characterization},
thereby establishing the uniqueness of
the associated epr-sequence (11b).

\textit{Case 2}:
$\sigma  \neq  1]0 10 \OL{10} 11 0 \OL{0}$.
Then, as $\sigma \neq 1]0110 \OL{0}$,
the epr-sequences (11b) and (11c) in
Theorem \ref{N in every 3 characterization}
are not associated with $\sigma$.
Also, it is clear that
$\sigma$ is not associated with the epr-sequence (11a) in
Theorem \ref{N in every 3 characterization}.
Since $\sigma \neq 0]110  \OL{0}$,
the epr-sequences (5a) and (5b) in
Theorem \ref{N in every 3 characterization}
are not associated with $\sigma$.
Thus far we have that
$\sigma$ is not associated with any of the epr-sequences
(5a), (5b), (11a), (11b) or (11c).
Hence, $\sigma$ must be one of the pr-sequences
(1)--(4), (6)--(10) or (12)--(15) in
Theorem \ref{No-111 characterization}.
Now, by considering all the possible cases,
one easily verifies that an attainable epr-sequence
associated with $\sigma$, which must be listed in
Theorem \ref{N in every 3 characterization},
is unique.
\epf

%%%\begin{cor}\label{Uniqueness: Nonsingular}
%%%Let $n \geq 3$ and $r_0]r_1 \cdots r_{n-1}1$ be a pr-sequence
%%%that is attainable by a real symmetric matrix.
%%%Suppose $r_0]r_1 \cdots r_{n-1}1$ does not contain three consecutive $1$s.
%%%Then there is a unique attainable epr-sequence associated with $r_0]r_1 \cdots r_{n-1}1$.
%%%\end{cor}
%%
%%
%%%Each of the pr-sequences
%%%$1]1101$ and $1]11011$,
%%%which have three consecutive $1$s,
%%%are associated with a unique attainable epr-sequence,
%%%namely $\tt{SSNA}$ and $\tt{SSNAA}$, respectively.
%%%%For $n \leq 5$ these are the only pr-sequences with
%%%%that property.

\section*{Acknowledgements}
$\null$
The author expresses his gratitude to Dr.\ Leslie Hogben,
for introducing him to this problem, and for her guidance.

%%%%%%%%%%%%%% References %%%%%%%%%%%%%%%%%%

\end{document}